\documentclass[english]{article}
\usepackage{babel}
\usepackage{array}
\usepackage{url}
\usepackage{multirow}
\usepackage{amsmath}
\usepackage{mismath}
\usepackage{amsthm}
\usepackage{graphicx}
\usepackage{float}
\usepackage{epstopdf}
\usepackage{epsfig}
\usepackage{tabularx,booktabs}
\newcolumntype{Y}{>{\centering\arraybackslash}X}
\usepackage[authoryear]{natbib}
\usepackage[unicode=true,pdfusetitle,
bookmarks=true,bookmarksnumbered=false,bookmarksopen=false,
breaklinks=false,pdfborder={0 0 1},backref=false,colorlinks=false]
{hyperref}
\hypersetup{
	colorlinks,urlcolor=blue,linkcolor=blue,anchorcolor=blue,citecolor=blue}
\usepackage{breakurl}

\makeatletter

\theoremstyle{plain}
\newtheorem{thm}{\protect\theoremname}
\theoremstyle{remark}
\newtheorem{rem}{\protect\remarkname}
\theoremstyle{plain}

\theoremstyle{definition}

\newcommand{\be}{\begin{equation}}
	\newcommand{\ee}{\end{equation}}
\newcommand{\bey}{\begin{eqnarray}}
	\newcommand{\eey}{\end{eqnarray}}
\newcommand{\beyn}{\begin{eqnarray*}}
	\newcommand{\eeyn}{\end{eqnarray*}}
\usepackage{color}
\usepackage{soul,xcolor}
\setstcolor{red}
\newif\ifApproveEdit
\ApproveEditfalse
\ifApproveEdit
  
  \newcommand{\DEL}[1]{\iffalse{#1}\fi}
\else
  
  \newcommand{\DEL}[1]{\st{#1}}
\fi

\usepackage{graphics}
\usepackage{epstopdf}

\renewcommand{\baselinestretch}{1.5} 

\usepackage{bm}

\usepackage[noblocks]{authblk}
\def\singlespace{\def\baselinestretch{1}\@normalsize}

\makeatother

\providecommand{\examplename}{Example}
\providecommand{\remarkname}{Remark}
\providecommand{\corollaryname}{Corollary}

\providecommand{\theoremname}{Theorem}

\global\long\def\b#1{{\bf \bm{\mathit{#1}}}}

\global\long\def\bzero{\b 0}

\global\long\def\ba{\b a}
\global\long\def\bA{\b A}

\global\long\def\bb{\b b}
\global\long\def\bB{\b B}

\global\long\def\bI{\b I}

\global\long\def\bJ{\b J}

\global\long\def\bn{\b n}

\global\long\def\bX{\b X}

\global\long\def\by{\b y}

\global\long\def\bz{\b z}

\global\long\def\bu{\b u}

\global\long\def\bw{\b w}

\global\long\def\bmu{\b{\mu}}

\global\long\def\bxi{\b{\xi}}

\global\long\def\bOmega{\b{\Omega}}

\global\long\def\bSigma{\b{\Sigma}}

\global\long\def\bOmega{\b{\Omega}}

\global\long\def\calK{\mathcal{K}}

\global\long\def\calN{\mathcal{N}}

\global\long\def\calL{\mathcal{L}}

\global\long\def\h#1{\hat{#1}}

\global\long\def\hd{\h d}

\global\long\def\halpha{\h{\alpha}}

\global\long\def\hbeta{\h{\beta}}

\global\long\def\hbOmega{\h{\bOmega}}

\global\long\def\t#1{\tilde{#1}}

\global\long\def\tT{\t T}

\global\long\def\barby{\bar{\by}}

\global\long\def\barbu{\bar{\bu}}

\global\long\def\N{\calN}

\global\long\def\T{\top}

\global\long\def\tr{\operatorname{tr}}

\global\long\def\diag{\operatorname{diag}}

\global\long\def\E{\operatorname{E}}

\global\long\def\Var{\operatorname{Var}}

\global\long\def\Cov{\operatorname{Cov}}

\global\long\def\ARE{\operatorname{ARE}}

\global\long\def\convP{\stackrel{P}{\longrightarrow}}

\global\long\def\convL{\stackrel{\calL}{\longrightarrow}}

\global\long\def\dequ{\stackrel{d}{=}}

\global\long\def\iidsim{\stackrel{\text{i.i.d.}}{\sim}}

\global\long\def\vs{\text{vs}}

\newcommand{\barbw}{\bar{\bw}}
\newcommand{\calU}{\mathcal{U}}
\newcommand{\calR}{\mathcal{R}}
\newcommand{\calB}{\mathcal{B}}
\newcommand{\sinf}{\mbox{Inf}}
\newcommand{\TS}{\mbox{\tiny S}}
\newcommand{\LC}{\mbox{\tiny LC}}
\newcommand{\CLX}{\mbox{\tiny CLX}}
\newcommand{\CZZW}{\mbox{\tiny CZZW}}
\newcommand{\NEW}{\mbox{\tiny NEW}}

\usepackage{xr}
\usepackage{soul,xcolor,cancel}
\usepackage{color}

\setstcolor{red}
\newif\ifApproveEdit
\ApproveEditfalse
\ifApproveEdit

\newcommand{\del}[1]{\iffalse{#1}\fi}
\else

\newcommand\del[2][red]{\setbox0=\hbox{$#2$}%
	\rlap{\raisebox{.45\ht0}{\textcolor{#1}{\rule{\wd0}{1pt}}}}#2}
\fi

\newif\ifApproveEdit
\ApproveEditfalse
\ifApproveEdit

\else

\fi

\begin{document}
	
	\global\long\def\GCBL{\operatorname{GCBL}}

	\author[a]{\small Jin-Ting Zhang}
	\author[a]{\small Jingyi Wang}
	\author[b]{\small Tianming Zhu}
	
	\affil[a]{\footnotesize Department of Statistics and Data Science, National University of Singapore, 
		Singapore}
	\affil[b]{\footnotesize Mathematics and Mathematics Education, National Institute of Education, Nanyang Technological University, 
		Singapore}
	
	\title{\Large Two-Sample Test for  High-Dimensional  Covariance Matrices: a normal-reference approach}
	\maketitle
	
	\vspace{-1cm}
	
	\begin{abstract}
	Testing the equality of the covariance matrices of two high-dimensional samples is a fundamental inference problem in statistics.
	Several tests have been proposed but they are either too liberal or too conservative when the required assumptions are not satisfied which attests that they are not always applicable in real data analysis. To overcome this difficulty, a normal-reference test is proposed and studied in this paper. It is shown that under some regularity conditions and the null hypothesis,  the proposed test statistic and a chi-square-type mixture have the same limiting distribution. It is then justified to approximate the null distribution of the proposed test statistic using that of the chi-square-type mixture. The distribution of the chi-square-type mixture can be well approximated using a three-cumulant matched chi-square-approximation with its approximation parameters consistently estimated from the data. The asymptotic power of the proposed test under a local alternative is also established. Simulation studies and a real data example demonstrate that in terms of size control, the proposed test outperforms the existing competitors substantially.
	\end{abstract}

\noindent{\bf KEY WORDS}: chi-square-type mixtures; high-dimensional data; three-cumulant matched  chi-square-approximation; equal covariance matrices.

\noindent{\bf Short Title}:  Test for  High-Dimensional  Covariance Matrices.
	\section{Introduction}\setcounter{equation}{0}\label{intro.sec}
	 Over the past few decades, with the availability of new technologies in data collection and storage, it has become common in economics, finance, and marketing to collect data on a large number of features for a limited number of individuals. In this situation,  the data dimension $p$ is close to or even much  larger than the total sample size $n$.  This ``large $p$, small $n$''  feature may cause most standard procedures
	to break down.  We refer to this kind of data as high-dimensional data and a problem about them  as a  ``large $p$, small $n$'' problem. An enduring recent interest in statistics is to
	compare the covariance  matrices of two high-dimensional populations since many statistical procedures rely on the fundamental assumption of equal covariance matrices. Examples include   \cite{bai1996effect}, \cite{zhangguozhoucheng2020JASA}, and \cite{zhang2019tsamsi} among others.   It is also common
	in  economics and finance.	The work of this paper is motivated by a financial data set provided by the Credit Research Initiative of National University of Singapore (NUS-CRI). 
    	The Probability of Default (PD) measures the likelihood of an obligor being unable to honor its financial obligations and is the core credit product of the NUS-CRI corporate default prediction system built on the forward intensity model of \cite{duan2012multiperiod}.  A major motivation for this study comes from the need to identify one of  the key distribution characteristics, that is, the covariance matrices, of daily PD which are significantly different with respect to emerging and developed financial markets. We are interested to check whether the daily PD in 2007--2008 of an emerging economy (mainland of China) and that of a developed economy (Hong Kong) have the same covariance matrix. The data set consists of $n_1=235$ companies listed in mainland of China and $n_2=153$ companies listed in Hong Kong, whereas each company has $p=522$ daily PD values in 2007--2008. Thus, it is a two-sample equal-covariance matrix testing problem for high-dimensional data.

	
	Mathematically, a two-sample  high-dimensional equal-covariance matrix testing  problem can be described as follows.  Suppose we have two independent high-dimensional samples:
	\begin{equation}\label{twosamp.sec1}
		\by_{i1},\ldots,\by_{in_i} \ \ \mbox{are i.i.d. with } \E(\by_{i1})=\bmu_i, \Cov(\by_{i1})=\bSigma_i,\; i=1,2,
	\end{equation}
	where the dimension $p$ is very large, and may
	be much larger than the total sample size $n=n_1+n_2$.  Of interest is to test whether the two covariance matrices  are equal:
	\begin{eqnarray}
		H_0 : \ \  \bSigma_1 = \bSigma_2 \ \ \vs\ \  H_1:\ \ \bSigma_1 \neq \bSigma_2.
		\label{covtest.sec2}
	\end{eqnarray}
	The above problem is not new and it has stimulated several works  in the literature in recent years.  Up to our knowledge,  \cite{Schott2007em2cov} may be the first author who studied  (\ref{covtest.sec2}) and proposed a test obtained via constructing an unbiased estimator to the usual squared Frobenius norm of the covariance  matrix difference $\bSigma_1-\bSigma_2$. Besides the  normality assumption,  \cite{Schott2007em2cov} imposed some other strong assumptions such that his test statistic is asymptotically normally distributed. However, in real data analysis, it is often difficult to check whether these strong assumptions are satisfied.  When these strong assumptions are not satisfied, we can show that \cite{Schott2007em2cov}'s  test statistic will not tend to normal asymptotically. In fact,  as shown in  Table~\ref{size1.tab} in Section~\ref{sim1.subsec},  \cite{Schott2007em2cov}'s test is very liberal: when the nominal size is 5\%,  for Gaussian data,  its empirical sizes can be as large as $10.96\% $ and for non-Gaussian data, its empirical sizes can be as large as 43.28\%.  It is not surprising that \cite{Schott2007em2cov}'s test performs much worse for non-Gaussian data than for Gaussian data since it  was proposed for Gaussian data only.  Without imposing the normality assumption,  \cite{LC2012} proposed a test via constructing an unbiased estimator to the usual squared Frobenius norm of the covariance  matrix difference $\bSigma_1-\bSigma_2$ using U-statistics. Under some strong assumptions,  \cite{LC2012} showed that their test statistic is also asymptotically normal. However, when these strong assumptions are not satisfied, \cite{LC2012}'s test is also very liberal as shown in  Table~\ref{size1.tab}. It is seen that the empirical sizes of \cite{LC2012}'s test can be as large as 13.61\%. Note that both \cite{Schott2007em2cov}
		and \cite{LC2012}'s tests are $L_2$-type, which may be less powerful when the entries of the covariance matrix difference $\bSigma_1-\bSigma_2$ are sparse. To overcome this difficulty,  \cite{cai2013} proposed an $L_{\infty}$-type test. They showed that under some regularity conditions,  their test statistic is asymptotically  an extreme-value distribution of Type I. Unfortunately, the simulation results in Table~\ref{size1.tab} indicate that \cite{cai2013}'s test is very conservative with its empirical sizes being  as small as $0.00\%$. This is a very undesired property. To overcome this drawback,  \cite{CZZW2017} proposed to approximate the finite-sample distribution of the $L_{\infty}$-type test using a bootstrap method. As indicated by Table~\ref{size1.tab}, in terms of size control,  the  bootstrap test of \cite{CZZW2017}  does work well  when the data are symmetrically distributed and are  moderately and highly correlated. However, it is still quite conservative with its empirical sizes being as small as 0.77\%  when the data are skewed and are nearly uncorrelated. From the above discussion, it is seen that the existing tests mentioned above cannot control their sizes well generally.
		
		To overcome the above size-control  problem, in this paper, we propose and study  a normal-reference test for (\ref{covtest.sec2}).  We first
		transform the two-sample equal-covariance matrix testing problem (\ref{covtest.sec2}) into a two-sample equal-mean vector testing  problem based on two induced high-dimensional samples obtained from the original two high-dimensional samples (\ref{twosamp.sec1}) with the help of  the well-known  Kronecker operator. To the best of our knowledge, this technique is new  for the two-sample high-dimensional equal-covariance matrix testing problem. Following \cite{chen2010two} and \cite{ZhangZhu2021HD2sBF3c}, a U-statistic based test statistic is constructed using the two induced high-dimensional samples. Under some regularity conditions and the null hypothesis, it is shown that the distribution difference between  the proposed test statistic and a chi-square-type mixture is   with order
		$\mathcal{O}(1/n_1+1/n_2)$ where $n_1$ and $n_2$ are the group sample sizes of the two samples. That is, the proposed test statistic and a chi-square-type mixture have the same limiting distribution. The limiting distribution  can be normal or non-normal, depending on how the components of the high-dimensional data are correlated. Therefore, it is not always applicable to approximate the null distribution of the proposed test statistic using a normal distribution as done in \cite{Schott2007em2cov} and \cite{LC2012} among others. On the other hand,  it is  justified that the null distribution of the proposed test can be well approximated by that of the chi-square-type mixture which is obtained when the two induced high-dimensional samples are normally distributed. We then term the proposed test as a normal-reference test. To approximate the distribution of the chi-square-type mixture, following \cite{ZhangZhu2021HD2sBF3c}, we apply the three-cumulant matched chi-square-approximation of \cite{zhang2005approximate} with the approximation parameters consistently estimated from the data. Under a local alternative, the asymptotic power of the proposed test is established. Two simulation studies and a real data example demonstrate that in terms of size control, the proposed test outperforms the existing competitors proposed by \cite{Schott2007em2cov}, \cite{LC2012}, \cite{cai2013}, and \cite{CZZW2017} substantially.
	
	The rest of the paper is organized as follows. The  main results are presented  in Section~\ref{main.sec}. Simulation studies and  an application to the financial data mentioned above  are given in Sections~\ref{simul.sec} and~\ref{appl.sec} respectively. Some concluding remarks are
	given in Section~\ref{remark.sec}. Technical proofs of the main results  are outlined in the Appendix.

\section{Main Results}\label{main.sec}

\subsection{Test statistic}

Without loss of generality and for simplicity, throughout this section,
we assume $\bmu_1=\bmu_2=\bzero$ since in this paper, we are concerned with  equal-covariance matrix testing problems only. In practice, it is often sufficient to replace $\by_{ij}, j=1,\ldots, n_i;\;i=1,2$  with $\by_{ij}-\barby_i, j=1,\ldots, n_i;\;i=1,2$  where $\barby_1$ and $\barby_2$  are  the usual group sample mean vectors of the two samples (\ref{twosamp.sec1}) when $\bmu_1$ and $\bmu_2$ are  actually  not equal to $\bzero$.
Under this assumption,  we can easily transform the  equal-covariance matrix  testing problem (\ref{covtest.sec2}) based on the two samples (\ref{twosamp.sec1})  into an  equal-mean vector testing problem based on two induced high-dimensional samples.

Let $vec(\bA)$ denote a column vector obtained via stacking  the column vectors of a matrix $\bA$ one by one. We have $vec(\by\by^{\top})=\by\otimes \by$ where $\otimes$ denotes the well-known  Kronecker operator  and $\by$ is a column vector. Then the equal-covariance matrix testing problem (\ref{covtest.sec2}) can be equivalently  written as the following equal-mean vector testing problem:
\be\label{Ncovtest.sec2}
H_0: vec(\bSigma_1)=vec(\bSigma_2)\quad \vs \quad H_1: vec(\bSigma_1)\neq vec(\bSigma_2),
\ee
based on the following two induced  samples
\be\label{newsamp.sec2}
\bw_{ij}=vec(\by_{ij}\by_{ij}^{\top})=\by_{ij}\otimes\by_{ij}, j=1,\ldots, n_i; \; i=1,2.
\ee
Following \citet{chen2010two} and \cite{ZhangZhu2021HD2sBF3c},  a U-statistic based  test statistic for (\ref{Ncovtest.sec2}) can be constructed  as
\[
T_{n,p}=\frac{2}{n_1(n_1-1)}\sum_{1\le i<j \le n_1} \bw_{1i}^{\top}\bw_{1j}+\frac{2}{n_2(n_2-1)}\sum_{1\le i<j \le n_2} \bw_{2i}^{\top}\bw_{2j}-\frac{2}{n_1n_2}\sum_{i=1}^{n_1} \sum_{j=1}^{n_2} \bw_{1i}^{\top}\bw_{2j}.
\]
It is easy to show that $\E(T_{n,p})=\|vec(\bSigma_1-\bSigma_2)\|^2=\tr\{(\bSigma_1-\bSigma_2)^2\}$, the usual squared Frobenius norm of the covariance  matrix difference $\bSigma_1-\bSigma_2$, where $\|\ba\|$ denotes the usual $L^2$-norm of a vector $\ba$. Hence it is justified to use $T_{n,p}$ to test (\ref{Ncovtest.sec2}). Let
\be\label{hOmega.sec2}
\barbw_i=n_i^{-1}\sum_{j=1}^{n_i} \bw_{ij}, \mbox{ and } \hbOmega_i=(n_i-1)^{-1}\sum_{j=1}^{n_i} (\bw_{ij}-\barbw_i)(\bw_{ij}-\barbw_i)^{\top}, i=1,2,
\ee
be the usual sample mean vectors and the usual sample covariance matrices of the two induced samples (\ref{newsamp.sec2}). By some simple algebra, we can equivalently write
\be\label{nstat.sec2}
T_{n,p}=\|\barbw_1-\barbw_2\|^2-\tr(\hbOmega_n),
\ee
where   $\hbOmega_n=\hbOmega_1/n_1+\hbOmega_2/n_2$.
Now in order to test (\ref{Ncovtest.sec2}), we need to derive the null distribution of $T_{n,p}$ defined in (\ref{nstat.sec2}). For this purpose, we set
\be\label{newsampA.sec2}
\bu_{ij}=\bw_{ij}-vec(\bSigma_i),\; j=1,\ldots, n_i; \;i=1,2,
\ee
and then we have $\E(\bu_{ij})=\bzero$ and $\Cov(\bu_{ij})=\bOmega_i,\;j=1,\ldots, n_i; \; i=1,2$.  It follows that  we can express $T_{n,p}$ as
\be\label{Tdecomp.sec2}
T_{n,p}=T_{n,p,0}+2S_{n,p}+\tr\{(\bSigma_1-\bSigma_2)^2\},
\ee
where
\be\label{Tstat0.sec2}
T_{n,p,0}=\|\barbu_1-\barbu_2\|^2-\tr(\hbOmega_n), \;\mbox{ and } \; S_{n,p}=(\barbu_1-\barbu_2)^{\T}\{vec(\bSigma_1-\bSigma_2)\},
\ee
with $\barbu_1$ and $\barbu_2$ are the usual sample mean vectors of the two centralized samples (\ref{newsampA.sec2}). It is easy to see that $T_{n,p,0}$ has the same distribution as that of $T_{n,p}$ under the null hypothesis.
Notice that we can write $T_{n,p,0}$ as
\[
T_{n,p,0}=\frac{2}{n_1(n_1-1)}\sum_{1\le i<j \le n_1} \bu_{1i}^{\top}\bu_{1j}+\frac{2}{n_2(n_2-1)}\sum_{1\le i<j \le n_2} \bu_{2i}^{\top}\bu_{2j}-\frac{2}{n_1n_2}\sum_{i=1}^{n_1} \sum_{j=1}^{n_2} \bu_{1i}^{\top}\bu_{2j},
\]
which is a quadratic form of the two centralized samples (\ref{newsampA.sec2}). For simplicity, set $\calU=\{\bu_{ij}=\bw_{ij}-vec(\bSigma_i), j=1,\ldots, n_i; \;i=1,2\}$, and we can write $T_{n,p,0}=Q(\calU)$. By some simple algebra, we have $\E(T_{n,p,0})=0$, and
\be\label{Tnsig2.sec2}
\sigma_T^2=\Var(T_{n,p,0})=2\left\{\frac{\tr(\bOmega_1^2)}{n_1(n_1-1)}+\frac{2\tr(\bOmega_1\bOmega_2)}{n_1n_2}+\frac{\tr(\bOmega_2^2)}{n_2(n_2-1)}\right\}.
\ee

In particular, when  the two induced  samples (\ref{newsamp.sec2}) are normally distributed, we have $T_{n,p,0}^*=Q(\calU^*)$, where $\calU^*=\{\bu_{ij}^*=\bw_{ij}^*-vec(\bSigma_i), j=1,\ldots, n_i;\; i=1,2\}$ denotes two centralized samples which are independent from $\calU$ but  we have $\bw_{ij}^*, j=1,\ldots, n_i \iidsim N_{p^2}(vec(\bSigma_i),\bOmega_i),\; i=1,2$.  Obviously, it follows that $\E(T_{n,p,0}^*)=0$ and $\Var(T_{n,p,0}^*)=\sigma_T^2$.  We call the distribution of $T_{n,p,0}^*$ as the normal-reference distribution of $T_{n,p,0}$. Throughout this paper, let $\chi_{v}^{2}$ denote a central chi-square distribution with $v$ degrees of freedom and  $\dequ$ denote equality in distribution.  For any given $n$ and $p$, it is easy to show  that $T_{n,p,0}^*$ has the same distribution as that of a chi-square-type mixture as follows:
\be\label{Tchisq.sec2}
T_{n,p,0}^*\dequ\sum_{r=1}^{p^2} \lambda_{n,p,r}A_r-\left\{\frac{\sum_{r=1}^{p^2}\lambda_{1r} B_{1r}}{n_1(n_1-1)}+\frac{\sum_{r=1}^{p^2}\lambda_{2r} B_{2r}}{n_2(n_2-1)}\right\},
\ee
where   $\lambda_{1r}, r=1,\ldots, p^2$ and $\lambda_{2r}, r=1,\ldots, p^2$ are the eigenvalues of $\bOmega_1$ and $\bOmega_2$, respectively, and
$A_{r}\iidsim\chi_1^2$, $B_{1r} \iidsim\chi_{n_1-1}^2$, and $B_{2r}\iidsim\chi_{n_2-1}^2$ are mutually independent.

\subsection{Asymptotic null distribution}

For further theoretical discussion, following \citet{wangxu2021}, we introduce a norm which measures the difference between two probability measures.  For two probability measures $\nu_1$ and $\nu_2$ on $\calR$, let $\nu_1-\nu_2$ denote  the signed measure such that for any Borel set $A$, $(\nu_1-\nu_2)(A) = \nu_1(A)-\nu_2(A)$. Let $\calB_b^3(\calR)$  denote
the class of bounded functions with continuous derivatives up to order $3$. It is known that a sequence of
random variables $\{x_i\}_{i=1}^{\infty}$ converges weakly to a random variable $x$ if and only if for every $f\in\calB_b^3(\calR)$, we have
$\E\{f(x_i)\}\to\E\{f(x)\}$; see \citet{wangxu2021} for some details. We use this property to give a definition of the weak convergence in $\calR$. For a function $f\in\calB_b^3(\calR)$, let $f^{(r)}$  denote the $r$th derivative of $f, r=1,2,3$. For a finite signed measure $\nu$ on $\calR$, we define the norm $\|\nu\|_3$  as $\sup_f\int_{\calR} f(x) \nu(dx)$
where the supremum is taken over all $f\in\calB_b^3(\calR)$ such that $\sup_{x\in\calR} |f^{(r)}(x)|\le 1, r=1,2,3$. It is straightforward
to verify that $\|\cdot\|_3$  is indeed a norm. Also, a sequence of probability measures $\{\nu_n\}_{n=1}^{\infty}$  converges weakly
to a probability measure $\nu$ if and only if $\|\nu_n-\nu\|_3\to 0$ as $n\to\infty$. Throughout this paper, for simplicity, we often denote $\{\E(\bX)\}^2$ and $\{\Var(\bX)\}^2$  as $\E^2(\bX)$ and $\Var^2(\bX)$, respectively. Let $\lambda_{n,p,r},r=1,\ldots, p^2 $ are  the eigenvalues of
\[
\bOmega_n=\Cov(\barbw_1-\barbw_2)=\bOmega_1/n_1+\bOmega_2/n_2,
\]
and set $c_{n,p,r}=\lambda_{n,p,r}/\sqrt{\tr(\bOmega_n^2)}$. For further study, we impose the following conditions:

\begin{description}
	\item[{C1.}] As $n\to\infty$, we have $n_1/n\to \tau\in (0,\infty)$.
	\item[{C2.}] There is a universal constant $1<\gamma<\infty$ such that for all $q\times p^2$ real matrix $\bB$, we have
	$
	\E\|\bB\bu_{ij}\|^4\le \gamma \E^2(\|\bB\bu_{ij}\|^2), \;  \mbox{ for all $j=1,\ldots, n_i;\; i=1,2$}.
	$
	\item [{C3.}] As $n,p\to \infty$, we have $c_{n,p,r}\to c_r\ge 0$ for all  $r=1,2,\ldots $ uniformly.
	\item[{C4.}] As $n,p\to\infty$, we have $p^2/n\to c\in(0,\infty)$.
	
\end{description}
Condition C1 is regular for any two-sample testing problem. It requires that the two sample sizes $n_1$ and $n_2$ tend to infinity proportionally. To give some insight about Condition C2, we list the following remarks.
\begin{rem} When $\bB$ is a row vector, e.g., $\bB=\bb^{\top}$, Condition C2 implies  that the kurtosis of $\bb^{\top}\bu_{ij}$ is  bounded by $\gamma$ for all $\bb$:
	$
	kurt(\bb^{\top}\bu_{ij})=\E(\bb^{\top}\bu_{ij})^4/\Var^2(\bb^{\top}\bu_{ij})\le \gamma.
	$
	That is,  the kurtosis of $\bu_{ij}=\bw_{ij}-vec(\bSigma_i)$ is uniformly  bounded  in any  projection direction for all $j=1,\ldots, n_i;\; i=1,2$. \cite{westfall2014kurtosis} clarified that the value of kurtosis is related to the tails of the distribution. 
	Therefore, Condition C2 essentially  requires that the distribution of $\bu_{ij}$ is not extremely  tailed in any projection direction. We shall see this is a fairly weak  condition.
\end{rem}

\begin{rem} We have
	$
	\E\|\bB\bu_{ij}\|^4=\Var\|\bB\bu_{ij}\|^2+\E^2(\|\bB\bu_{ij}\|^2),
	$
	which explains why $\gamma$ must be larger than $1$. Together with  Condition C2,  imply  that the variances  of $\|\bB\bu_{ij}\|^2$'s  are  uniformly bounded by their   squared means $\E^2(\|\bB\bu_{ij}\|^2)$ and the noise-to-signal ratios $(\Var\|\bB\bu_{ij}\|^2)^{1/2}/\E\|\bB\bu_{ij}\|^2$ are also uniformly bounded.
\end{rem}

\begin{rem} Condition C2 is automatically satisfied by any Gaussian samples with $\gamma=3$. Since $\Cov(\bu_{ij})=\bOmega_i$, we can express the   two centralized samples (\ref{newsampA.sec2}) as  $\bu_{ij}=\bOmega_i^{1/2}\bz_{ij}, j=1,\ldots, n_i; \;i=1,2$, where $\bz_{i1},\ldots, \bz_{in_i} $ are i.i.d.  with $\E(\bz_{ij})=\bzero$ and $\Cov(\bz_{ij})=\bI_{p^2}$. If the components $z_{ij1},\ldots, z_{ijp^2}$  of $\bz_{ij}$ are independent with $\E(z_{ijr}^4)=\Var(z_{ijr}^2)+1\le \Delta<\infty$ for all $j=1,\ldots, n_i; \;i=1,2;\; r=1,\ldots, p^2$, then it is easy to  show that Condition C2 is satisfied with $\gamma=\Delta$.
\end{rem}
Condition C3 ensures the existence of the limits of $c_{n,p,r}$ which are the eigenvalues of $\bOmega_n/\sqrt{\tr(\bOmega_n^2)}$. It is used to get the limiting distributions of the standardized versions of $T_{n,p,0}$ and $T_{n,p,0}^*$, namely,
\begin{equation*}\label{stdTstat.sec2}
	\tT_{n,p,0}=T_{n,p,0}/\sigma_T, \mbox{ and } \tT_{n,p,0}^*=T_{n,p,0}^*/\sigma_T,
\end{equation*}
where $\tT_{n,p,0}$ and $\tT_{n,p,0}^*$ have zero mean and unit variance. Condition C4 is imposed for studying the ratio-consistency of the estimators used in the proposed normal-reference test. Throughout this paper,  let $\calL(y)$  denote the distribution of a random variable $y$ and $\convL$ denote convergence in distribution.
We then have the following useful theorem.

\begin{thm} \label{Tnexp.thm} Under Condition C2, we have
	\[
	\|\calL(\tT_{n,p,0})-\calL(\tT_{n,p,0}^*)\|_3\le \frac{(2\gamma)^{3/2}}{3^{1/4}}\left(\frac{1}{n_1}+\frac{1}{n_2}\right)^{1/2},
	\]
	where $\gamma$ is defined in Condition C2.
\end{thm}

Theorem~\ref{Tnexp.thm} indicates that the distance between the distribution of $\tT_{n,p,0}$ and the distribution of $\tT_{n,p,0}^*$  is $\mathcal{O}\{(1/n_1+1/n_2)^{1/2}\}$, showing that the distributions of $\t{T}_{n,p,0}$ and $\t{T}_{n,p,0}^*$ are equivalent asymptotically.  Thus, Theorem~\ref{Tnexp.thm} actually provides a systematic theoretical justification for us to use the distribution $T_{n,p,0}^*$ to approximate the distribution of $T_{n,p,0}$.

\begin{thm} \label{Tnlim.thm} Under Conditions C1--C3, as $n,p\rightarrow \infty$, we have $\tT_{n,p,0}^*\convL \zeta$ with
	\be\label{Tnlim.sec2}
	\zeta\dequ (1-\sum_{r=1}^{\infty}c_r^2)^{1/2}z_0+2^{-1/2}\sum_{r=1}^{\infty} c_r(z_r^2-1),
	\ee
	where $z_0,z_1,z_2,\ldots$ are i.i.d. $N(0,1)$, and $c_r,r=1,2,\ldots,$ are defined in Condition C3.
\end{thm}
Theorem \ref{Tnlim.thm} gives a unified expression  of the possible asymptotic distributions of $\tT_{n,p,0}^*$,  a weighted sum of a standard normal random variable and a sequence of  centered chi-square  random variables. From Fatou's Lemma and Condition C3, $\sum_{r=1}^{\infty}c_r^2 \leq \lim_{n,p \to \infty}\sum_{r=1}^{p}c_{n,p,r}^2=1$. This shows that  $\sum_{r=1}^{\infty}c_r^2 \in [0,1]$. Some remarks are given below to reveal some special cases of the possible  distributions of $\zeta$ (\ref{Tnlim.sec2}).

\begin{rem}\label{nordist.rem} We have $\zeta\dequ z_0\sim N(0,1)$ when $\sum_{r=1}^{\infty}c_r^2 = 0 $,  equivalently,  $c_r=0, r=1,2,\ldots$ which
		holds when one of the following conditions holds: as $n,p\to\infty$,
		\be\label{norconds.sec2}
		\begin{array}{ll}
			\lambda_{n,p,max}^2 = \lito\{\tr(\bOmega_n^2)\}, &\mbox{\cite{bai1996effect}}, \\
			\tr(\bOmega_n^4) = \lito\{\tr^2(\bOmega_n^2)\},   &\mbox{\cite{chen2010two}},\\
			\tr(\bOmega_n^\ell)/p^2 \to a_\ell  \in (0,\infty), \ell=1,2,3,  &\mbox{\cite{srivastava2008test}},
		\end{array}
		\ee
		where $\lambda_{n,p,max}$ denotes the largest eigenvalue of $\bOmega_n$.
\end{rem}

\begin{rem} We have $\zeta\dequ \sum_{r=1}^{\infty} c_r(z_r^2-1)$, a weighted sum of centered chi-square random variables when $\sum_{r=1}^{\infty}c_r^2 = 1$, which holds under Condition C3 and when the limit and summation operations in $\lim_{n,p\to\infty} \sum_{r=1}^{p^2} c_{n,p,r}^2$ are exchangeable, that is, when
 $\lim_{n,p\to\infty} \sum_{r=1}^{p^2} c_{n,p,r}^2=\sum_{r=1}^{\infty}\lim_{n,p\to\infty} c_{n,p,r}^2.$
\end{rem}

\begin{rem} The above two remarks  indicate that the null limiting distribution of $T_{n,p}$ can be normal or non-normal.
		However, in practice, it is often challenging to check whether $\sum_{r=1}^{\infty}c_r^2=0$ or $\sum_{r=1}^{\infty}c_r^2=1$. Therefore, it is not always appropriate to  use the normal approximation to the null distribution of $T_{n,p}$ as done by \cite{Schott2007em2cov} and \cite{LC2012}. Roughly speaking, when the $p^2$ components of the two induced  samples (\ref{newsamp.sec2}) are nearly uncorrelated, we have $\sum_{r=1}^{\infty}c_r^2=0$ and when they  are moderately or highly correlated, we have $\sum_{r=1}^{\infty}c_r^2=1$.   Fortunately, in the proposed normal-reference test, we  approximate the null distribution of $T_{n,p}$ using the
		distribution of $T_{n,p,0}^*$, justified by Theorem~\ref{Tnexp.thm},   rather than using  its limiting distribution obtained in Theorem~\ref{Tnlim.thm}.  This is an advantage of the proposed  normal-reference test.
\end{rem}

\subsection{Null distribution approximation \label{implement.sec}}

Theorem~\ref{Tnexp.thm} shows that it is justified to approximate the distribution of $T_{n,p,0}$ using that of $T_{n,p,0}^*$. Notice from (\ref{Tchisq.sec2})
that $T_{n,p,0}^*$ is  generally skewed although under certain regularity  conditions, it can be asymptotically normally distributed as mentioned in Remark~\ref{nordist.rem}.
Since $T_{n,p,0}^*$  is  a $\chi^2$-type mixture with  positive and negative unknown coefficients, its distribution should not be approximated
using the Welch--Satterthwaite $\chi^2$-approximation as in \citet{zhangguozhouzhu2020}; rather we should  approximate its distribution using the 3-c matched
$\chi^2$-approximation (\citealt{zhang2005approximate}, \citealt{zhang2013analysis}). The key idea of  the 3-c matched $\chi^2$-approximation
is to approximate the distribution of $T_{n,p,0}^*$ using that of the following random variable
\[
R\dequ \beta_0+\beta_1\chi_d^2,
\]
where $\beta_0, \beta_1$, and $d$ are the approximation parameters with $d$ being called the approximate degrees of freedom of the 3-c matched $\chi^2$-approximation. They are determined
via matching the first three cumulants of $T_{n,p,0}^*$ and $R$. The first three cumulants of a random variable are its mean, variance, and third central moment (\citealt{zhang2005approximate}).
The first three cumulants of $R$ are given by  $\calK_1(R)=\beta_0+\beta_1 d$, $\calK_2(R)=2\beta_1^2 d$, and $\calK_3(R)=8\beta_1^3 d$ while the first three cumulants of $T_{n,p,0}^*$ are given  by $\calK_1(T_{n,p,0}^*)=\E(T_{n,p,0}^*)=0$,
\be\label{cumulants.sec2}
\begin{split}
	\calK_2(T_{n,p,0}^*)&=\Var(T_{n,p,0}^*)=\sigma_T^2=2\left\{\tr(\bOmega_n^2)+\frac{\tr(\bOmega_1^2)}{n_1^2(n_1-1)}+\frac{\tr(\bOmega_2^2)}{n_2^2(n_2-1)}\right\},\;\mbox{ and }\\
	\calK_3(T_{n,p,0}^*)&=\E({T_{n,p,0}^*}^3)=8\left\{\tr(\bOmega_n^3)-\frac{\tr(\bOmega_1^3)}{n_1^3(n_1-1)^2}-\frac{\tr(\bOmega_2^3)}{n_2^3(n_2-1)^2}\right\}.  
\end{split}
\ee
Matching  the first three cumulants of $T_{n,p,0}^*$ and $R$ then leads to
\be\label{betadf.sec2}
\beta_0=-\frac{2\calK_2^2(T_{n,p,0}^*)}{\calK_3(T_{n,p,0}^*)},\;\;
\beta_1=\frac{\calK_3(T_{n,p,0}^*)}{4\calK_2(T_{n,p,0}^*)},\;\;\mbox{and}\;\;
d=\frac{8\calK_2^3(T_{n,p,0}^*)}{\calK_3^2(T_{n,p,0}^*)}.
\ee
From (\ref{cumulants.sec2}), by some simple algebra, we have
\be\label{cumulant2.sec2}
\begin{split}
	\calK_2(T_{n,p,0}^*)&=2\left\{\frac{\tr(\bOmega_1^2)}{n_1(n_1-1)}+\frac{2\tr(\bOmega_1\bOmega_2)}{n_1n_2}+\frac{\tr(\bOmega_2^2)}{n_2(n_2-1)}\right\},\mbox{ and }\\
	\calK_3(T_{n,p,0}^*)&=8\left\{\frac{(n_1-2)\tr(\bOmega_1^3)}{n_1^2(n_1-1)^2}+\frac{3\tr(\bOmega_1^2\bOmega_2)}{n_1^2n_2}
	+\frac{3\tr(\bOmega_1\bOmega_2^2)}{n_1n_2^2}+\frac{(n_2-2)\tr(\bOmega_2^3)}{n_2^2(n_2-1)^2}\right\}.
\end{split}
\ee
It is seen  that  $\calK_2(T_{n,p,0}^*)>0$ and $\calK_3(T_{n,p,0}^*)>0$ since we always should have
$n_1, n_2>2$ and $\bOmega_i,i=1,2$ to be nonnegative. Thus, we always have $\beta_0<0, \beta_1>0$, and $d>0$. The negative value of $\beta_0$  is expected  since $T_{n,p,0}^*$ is a chi-square-type mixture with both positive and negative coefficients.  Note that the skewness of $T_{n,p,0}^*$ can be expressed as
\be\label{skewness.sec2}
\frac{\E({T_{n,p,0}^*}^3)}{\Var^{3/2}(T_{n,p,0}^*)}=\frac{\calK_3(T_{n,p,0}^*)}{\calK_2^{3/2}(T_{n,p,0}^*)}=\sqrt{8/d}.
\ee

\begin{rem}\label{abetadf.rem} For large $n_1$ and $n_2$, by (\ref{cumulants.sec2}), we have
	\[
	\calK_2(T_{n,p,0}^*)=2\tr(\bOmega_n^2)\{1+\lito(1)\}, \;\mbox{  and }\;\;  \calK_3(T_{n,p,0}^*)=8\tr(\bOmega_n^3)\{1+\lito(1)\}.
	\]
	Then by (\ref{betadf.sec2}), we have
	\[
	\beta_0=-\frac{\tr^2(\bOmega_n^2)}{2\tr(\bOmega_n^3)}\{1+\lito(1)\},\;\;
	\beta_1=\frac{\tr(\bOmega_n^3)}{\tr(\bOmega_n^2)}\{1+\lito(1)\},\;\;\mbox{and}\;\;
	d=\frac{\tr^3(\bOmega_n^2)}{\tr^2(\bOmega_n^3)}\{1+\lito(1)\}.
	\]
\end{rem}

	\begin{rem}
		The expression (\ref{skewness.sec2}) indicates that the value of $d$ can be used to quantify the skewness of $T_{n,p,0}^*$ and  when $\tT_{n,p,0}^*$ is asymptotically normal, $d$ must tend to $\infty$ and when $d$ is bounded, $\tT_{n,p,0}^*$ will not tend to normal. By the proof of Theorem 3 of \cite{zhangguozhoucheng2020JASA}, we have
			\be\label{dstar.sec2}
			d^{*-1}\le \frac{\tr(\bOmega_n^4)}{\tr^2(\bOmega_n^2)}\le \frac{\lambda_{n,p,\max}^2}{\tr(\bOmega_n^2)}\le d^{*-1/3},
			\ee
			where $d^*=\tr^3(\bOmega_n^2)/\tr^3(\bOmega_n^2)$. By Remark~\ref{abetadf.rem}, when $n_1$ and $n_2$ are large, we have $d=d^*\{1+\lito(1)\}$. Therefore, when $d\to\infty$, we have $d^*\to\infty$.  Then
			by (\ref{dstar.sec2}), it follows that  the first two conditions of (\ref{norconds.sec2}) hold. Therefore, we have $\tT_{n,p,0}^*$ tends to $N(0,1)$. Thus,
			$\tT_{n,p,0}^*$ is asymptotically normal if and only if $d\to \infty$, and hence we can use the value of $d$ to assess if the normal approximation to the null distribution of $T_{n,p}$ is adequate. This may also be regarded as an advantage for using the proposed normal-reference test with the 3-c matched $\chi^2$-approximation.
	\end{rem}
	To apply   the 3-c matched $\chi^2$-approximation,  we  need to estimate $\calK_2(T_{n,p,0}^*)$ and $\calK_3(T_{n,p,0}^*)$ consistently.
	Recall that the usual unbiased estimators of $\bOmega_1$ and $\bOmega_2$ are given by   $\hbOmega_1$ and $\hbOmega_2$  as  in  (\ref{hOmega.sec2}).  We first find an unbiased and ratio-consistent estimator of $\calK_2(T_{n,p,0}^*)$.  According to (\ref{cumulant2.sec2}), to obtain an unbiased and ratio-consistent estimator of $\calK_2(T_{n,p,0}^*)$, we need the  unbiased and  ratio-consistent estimators of $\tr(\bOmega_1^2), \tr(\bOmega_2^2)$, and $\tr(\bOmega_1\bOmega_2)$, respectively. By Lemma S.3 of  \citet{zhangguozhoucheng2020JASA}, the unbiased and ratio-consistent estimators of $\tr(\bOmega_i^2), i=1,2$ are given by
	\[
	\widehat{\tr(\bOmega_i^2)}=\frac{(n_i-1)^2}{(n_i-2)(n_i+1)}\left\{\tr(\hbOmega_i^2)-\frac{1}{n_i-1}\tr^2(\hbOmega_i)\right\}, i=1,2.
	\]
	By the proof of Theorem 2 of \citet{zhangguozhouzhu2020}, the unbiased and ratio-consistent estimator of $\tr(\bOmega_1\bOmega_2)$ is given by $\tr(\hbOmega_1\hbOmega_2)$. Therefore,
	based on (\ref{cumulant2.sec2}), the unbiased and ratio-consistent estimator of  $\calK_2(T_{n,p,0}^*)$ is given by
	\[
	\widehat{\calK_2(T_{n,p,0}^*)}=2\left\{ \frac{\widehat{\tr(\bOmega_1^2)}}{n_1(n_1-1)}+\frac{2\tr(\hbOmega_1\hbOmega_2)}{n_1 n_2}
	+\frac{\widehat{\tr(\bOmega_2^2)}}{n_2(n_2-1)}\right\}.
	\]
	We now  find an unbiased and ratio-consistent estimator of $\calK_3(T_{n,p,0}^*)$. According to (\ref{cumulant2.sec2}), to obtain an unbiased and ratio-consistent estimator of $\calK_3(T_{n,p,0}^*)$, we need the  unbiased and  ratio-consistent estimators of $\tr(\bOmega_1^3), \tr(\bOmega_2^3)$,  $\tr(\bOmega_1^2\bOmega_2)$, and $\tr(\bOmega_1\bOmega_2^2)$,  respectively.
	By Lemma 1 of \citet{zhangzhouguo2020onesamp}, under Condition C4 and when the two induced samples (\ref{newsamp.sec2}) are normally distributed,  the unbiased and ratio-consistent estimators of $\tr(\bOmega_1^3)$ and $\tr(\bOmega_2^3)$ are given by
	\[
	\widehat{\tr(\bOmega_i^3)}=\frac{(n_i-1)^4}{(n_i^2+n_i-6)(n_i^2-2n_i-3)}\left\{\tr(\hbOmega_i^3)-\frac{3\tr(\hbOmega_i)\tr(\hbOmega_i^2)}{(n_i-1)}
	+\frac{2\tr^3(\hbOmega_i)}{(n_i-1)^2}\right\},\;   i=1,2.
	\]
	By Lemma 1 of \cite{ZhangZhu2021HD2sBF3c},  when the two induced samples (\ref{newsamp.sec2}) are normally distributed,  the unbiased estimators of  $\tr(\bOmega_1^2\bOmega_2)$ and $\tr(\bOmega_1\bOmega_2^2)$ are given by
	\[
	\begin{split}
		\widehat{\tr(\bOmega_1^2\bOmega_2)}&=\frac{(n_1-1)}{(n_1-2)(n_1+1)}\left\{(n_1-1)\tr(\hbOmega_1^2\hbOmega_2)-\tr(\hbOmega_1\hbOmega_2)\tr(\hbOmega_1)\right\},\;\mbox{ and }\\
		\widehat{\tr(\bOmega_1\bOmega_2^2)}&=\frac{(n_2-1)}{(n_2-2)(n_2+1)}\left\{(n_2-1)\tr(\hbOmega_1\hbOmega_2^2)-\tr(\hbOmega_1\hbOmega_2)\tr(\hbOmega_2)\right\},
	\end{split}
	\]
	respectively. Under some regularity conditions and when the two induced samples (\ref{newsamp.sec2}) are normally distributed, \citet{Hyodo2020statprob} showed that the above estimators are also ratio-consistent for
	$\tr(\bOmega_1^2\bOmega_2)$ and $\tr(\bOmega_1\bOmega_2^2)$. Then  the unbiased and ratio-consistent estimator of $\calK_3(T_{n,p,0}^*)$ is given by
	\[
	\widehat{\calK_3(T_{n,p,0}^*)}= 8\left\{\frac{(n_1-2)\widehat{\tr(\bOmega_1^3)}}{n_1^2(n_1-1)^2}+\frac{3\widehat{\tr(\bOmega_1^2\bOmega_2)}}{n_1^2 n_2}+\frac{3\widehat{\tr(\bOmega_1\bOmega_2^2)}}{n_1 n_2^2}
	+\frac{(n_2-2)\widehat{\tr(\bOmega_2^3)}}{n_2^2(n_2-1)^2}\right\}.
	\]
	It follows that the ratio-consistent  estimators of $\beta_0, \beta_1$, and $d$ are  given by
	\be\label{hbetadf.sec2}
	\hbeta_0=-\frac{2\{\widehat{\calK_2(T_{n,p,0}^*)}\}^2}{\widehat{\calK_3(T_{n,p,0}^*)}},\;\;
	\hbeta_1=\frac{\widehat{\calK_3(T_{n,p,0}^*)}}{4\widehat{\calK_2(T_{n,p,0}^*)}},\;\;\mbox{and}\;\;
	\hd=\frac{8\{\widehat{\calK_2(T_{n,p,0}^*)}\}^3}{\{\widehat{\calK_3(T_{n,p,0}^*)}\}^2}.
	\ee
	For any nominal significance level $\alpha>0$, let
	$\chi_{d}^2(\alpha)$ denote the upper $100\alpha$ percentile of
	$\chi_{d}^2$. Then using  (\ref{hbetadf.sec2}),  the new normal-reference  test for  the two-sample equal-covariance matrix testing  problem
	(\ref{covtest.sec2}) is  then conducted via using  the approximate critical
	value $\hbeta_0+\hbeta_1\chi_{\hd}^2(\alpha)$ or  the approximate $p$-value
	$\Pr\{\chi_{\hd}^2\ge (T_{n,p}-\hbeta_0)/\hbeta_1\}$.
	
	In practice, one may often use the following  normalized version of $T_{n,p}$:
	\[
		\tT_{n,p}=\frac{T_{n,p}}{\sqrt{\widehat{\calK_2(T_{n,p,0}^*)}}}.
		\]
	Then to approximate the null distribution of  $T_{n,p}$ using the distribution  of $\hbeta_0+\hbeta_1\chi_{\hd}^2$ is equivalent to approximate the null distribution of $\tT_{n,p}$ using the distribution  of
	$(\chi_{\hd}^2-\hd)/\sqrt{2\hd}$. In this case, the new normal-reference   test using $\tT_{n,p}$ is then conducted via
	using the approximate critical value $\{\chi_{\hd}^{2}(\alpha)-\hd\}/\sqrt{2\hd}$
	or the approximate $p$-value $\Pr( \chi_{\hd}^{2}\ge \hd+\sqrt{2\hd}\tT_{n,p})$.

	\begin{rem}\label{bound.rem}  \citet{zhang2005approximate} showed that under some regularity conditions, the upper density approximation error bounds for approximating the distribution of $T_{n,p,0}^*$ using the normal approximation and the 3-c matched $\chi^2$-approximation  are
		$\mathcal{O}(d^{-1/2})$ and $\mathcal{O}(M)+\mathcal{O}(d^{-1})$ respectively  where $d$ is defined in (\ref{betadf.sec2}) and  $M=\frac{\tr(\bOmega_n^4)}{\tr^2(\bOmega_n^2)}\{1+\lito(1)\}$  as $n, p\to\infty$. Thus, for approximating the distribution of $T_{n,p,0}^*$, it is theoretically justified that the proposed test with the 3-c matched $\chi^2$-approximation is much more accurate than the tests proposed by \cite{Schott2007em2cov} and \cite{LC2012}  with the normal approximation.
	\end{rem}	

\subsection{Asymptotic power \label{power.sec}}

We now  consider the asymptotic  power of $T_{n,p}$  under the following local alternative:
\begin{equation}
	\Var(S_{n,p})=vec(\bSigma_1-\bSigma_2)^\top\bOmega_n vec(\bSigma_1-\bSigma_2)=\lito\{\tr(\bOmega_n^2)\}\quad \mbox{as \ensuremath{n,p\rightarrow\infty}},\label{powercond1.sec2}
\end{equation}
where $S_{n,p}$ is defined in (\ref{Tstat0.sec2}). Under Condition C1, as $n\to\infty$, we have
\be\label{Omega.sec2}
\bOmega_n=\{n\tau(1-\tau)\}^{-1}\bOmega\{1+\lito(1)\}, \mbox{ where } \bOmega=(1-\tau)\bOmega_1+\tau\bOmega_2.
\ee

\begin{thm} \label{Tnpower.thm} Assume that as $n, p\to\infty$,  $\hbeta_0,\hbeta_1$, and $\hd$ are ratio-consistent for $\beta_0,\beta_1$, and $d$, respectively.	Under Conditions C1--C4,  and the local alternative (\ref{powercond1.sec2}),  as $n, p\to\infty$, we have
	\be\label{powerA.sec2}
	\Pr\left\{T_{n,p}>\hbeta_0+\hbeta_1\chi_{\hd}^2(\alpha)\right\}
	=\Pr\left[\zeta \ge \frac{\chi_{d}^2(\alpha)-d}{\sqrt{2d}}-\frac{n\tau(1-\tau)}{\sqrt{2\tr(\bOmega^2)}}\tr\left\{(\bSigma_1-\bSigma_2)^2\right\}\right]\{1+\lito(1)\},
	\ee
	where $\zeta$ is defined in Theorem~\ref{Tnlim.thm}. In addition, as $d\to \infty$, the above expression can be further expressed as
	\be\label{powerB.sec2}
	\Pr\left\{T_{n,p}>\hbeta_0+\hbeta_1\chi_{\hd}^2(\alpha)\right\}=\Phi\left[-z_{\alpha}+\frac{n\tau(1-\tau)}{\sqrt{2\tr(\bOmega^2)}}\tr\left\{(\bSigma_1-\bSigma_2)^2 \right\}\right]\{1+\lito(1)\},
	\ee
	where  $z_{\alpha}$ denotes the upper $100\alpha$-percentile of $N(0,1)$.
\end{thm}
\section{Simulation Studies}\label{simul.sec}

In this section, we conduct two simulation studies to compare the performance of the proposed normal-reference test, denoted as $T_{\NEW}$,  against  several existing competitors for the two-sample high-dimensional equal-covariance matrix testing problem (\ref{covtest.sec2}), including the tests by \cite{Schott2007em2cov},  \cite{LC2012}, \cite{cai2013} and \cite{CZZW2017}, denoted as $T_{\TS}$, $T_{\LC}$, $T_{\CLX}$ and $T_{\CZZW}$, respectively. These existing competitors have been briefly reviewed in Section~\ref{intro.sec}.  Throughout this section, the nominal size $\alpha$ and the number of simulation runs are set as $5\%$ and $10,000$, respectively. The proportion of the number of rejections out of $10,000$ simulation runs for a test is known  as the empirical size or power of the test.  To assess  the overall performance of a test in maintaining the nominal size, we adopt the  average relative error (ARE) advocated  by \cite{Zhang2011tech}. It is calculated as $\ARE=100M^{-1}\sum_{j=1}^M|\halpha_j-\alpha|/\alpha$, where $\halpha_j,j=1,\ldots,M$ denotes the empirical sizes under $M$ simulation settings. A smaller ARE value of a test  indicates   a better overall performance of that test  in terms of size control.

\subsection{Simulation 1}\label{sim1.subsec}

In this simulation study,   we generate the two high-dimensional samples (\ref{twosamp.sec1}) using
$
\by_{ij} =\bmu+ \bSigma_i^{1/2}\bz_{ij},j=1,...,n_i;\;i=1,2,
$
where $\bz_{ij}=(z_{ij1},\ldots,z_{ijp})^\T, j=1,\ldots, n_i;\; i=1,2$ are i.i.d.  random variables with $\E(\bz_{ij})=\bzero$ and $\Cov(\bz_{ij})=\bI_p$. The $p$ entries of $\bz_{ij}$ are generated using  the following three models:

\begin{description}
	\label{sam.des}
	\item[Model 1:] $z_{ijk},k=1,\ldots,p \iidsim \N(0,1)$.
	\item[Model 2:] $z_{ijk}= u_{ijk}/\sqrt{5/3}$, with $u_{ijk},k=1,\ldots,p \iidsim t_5$.
	\item[Model 3:] $z_{ijk}= (u_{ijk}-1)/\sqrt{2}$, with $u_{ijk},k=1,\ldots,p \iidsim \chi^2_1$.
\end{description}

The above three models generate  $z_{ijk}$'s with three distributions: normal; nonnormal but symmetric; and nonnormal and skewed. Without loss of generality, we specify $\bmu = \bzero$. The covariance matrices $\bSigma_1$ and $\bSigma_2$ are specified as $\bSigma_i=\diag(\sigma_{i1}^2,\ldots,\sigma_{ip}^2)\{(1-\rho_i)\bI_p+\rho_i\bJ_p\},i=1,2$, where $\bJ_p$ is the $p\times p$ matrix of ones. The covariance matrix difference $\bSigma_1-\bSigma_2$ is determined by $\sigma_{ik}^2,k=1\ldots,p;\; i=1,2$ and $\rho_i,i=1,2$. In particular, $\sigma_{ik}^2$'s control the variances of the generated two samples (\ref{twosamp.sec1}) while $\rho_1$ and $\rho_2$ control their correlations. When $\sigma_{1k}^2=\sigma_{2k}^2,k=1,\ldots,p$ and $\rho_1=\rho_2$, we have $\bSigma_1=\bSigma_2$ so that the two samples   are generated under the null hypothesis in (\ref{covtest.sec2}). For simplicity, we set $\sigma_{1k}^2=\sigma_{2k}^2=4,k=1,\ldots,p$ and specify $\rho_1=\rho_2=\rho$ with  $\rho=0.25,0.5,$ and $0.9$ to compute the empirical sizes of the tests under consideration. To compare the powers of the tests, we set $\sigma_{1k}^2=4,k=1,\ldots,p$ and $\sigma_{2k}^2=3.5+u_k,k=1,\ldots,p$ where $u_k,k=1,\ldots,p$ are randomly generated from the standard uniform distribution $U[0,1]$. In addition, we set $\rho_1=0.5$ and  consider three cases of $\rho_2=\rho=0.25,0.5$, and $0.9$. Finally, we  set $p=50,100$, and $500$, and $\bn=(n_1,n_2)=\bn_1, \bn_2$, and $\bn_3$  with  $\bn_1=(50,80)$, $\bn_2=(80,120)$, and $\bn_3=(200,300)$, respectively.

\begin{table}
	\caption{\label{size1.tab}	Empirical sizes (in $\%$) of $T_{\TS}, T_{\LC}, T_{\CLX}, T_{\CZZW}$, and $T_{\NEW}$ in  Simulation 1 with $\rho_1=\rho_2=\rho$.}
	\centering
		\setlength{\tabcolsep}{3pt}	
		\scalebox{0.75}{
			\begin{tabularx}{1.3\textwidth}{@{\extracolsep{\fill}}cccccccccccccccccc}
				\hline
				\textbf{}& \textbf{}& \textbf{}& \multicolumn{5}{c}{$\rho = 0.25$}& \multicolumn{5}{c}{$\rho = 0.5$}& \multicolumn{5}{c}{$\rho = 0.9$}\\ \hline
				Model & $p$ & $\bn$ & $T_{\TS}$ &$T_{\LC}$ & $T_{\CLX}$ & $T_{\CZZW}$ & $T_{\NEW}$ & $T_{\TS}$ & $T_{\LC}$ & $T_{\CLX}$ & $T_{\CZZW}$ & $T_{\NEW}$ & $T_{\TS}$ & $T_{\LC}$ & $T_{\CLX}$ & $T_{\CZZW}$ & $T_{\NEW}$ \\ \hline
				\multirow{9}{*}{1} & \multirow{3}{*}{50} & $\bn_1$& 6.57 & 7.13 & 5.58 & 7.10 & 4.62 & 8.34 & 10.31 & 4.15 & 8.66 & 6.44 & 8.92 & 11.97 & 0.41 & 7.22 & 5.74\\
				&  & $\bn_2$ & 6.23 & 7.98 & 5.95 & 7.02 & 4.66 & 8.57 & 10.52 & 3.44 & 8.08 & 5.62 & 8.31 & 9.02 & 0.20 & 5.02 & 4.92\\
				&  & $\bn_3$& 8.09 & 9.78 & 4.44 & 5.06 & 4.48 & 7.83 & 9.44 & 2.59 & 4.63 & 4.74 & 10.35 & 11.98 & 0.21 & 5.90 & 6.51\\ \cline{2-18}
				& \multirow{3}{*}{100}& $\bn_1$& 7.14 & 8.45 & 6.21 & 8.55 & 5.41 & 6.65 & 9.00 & 4.31 & 10.49 & 4.28 & 8.11 & 9.75 & 0.32 & 6.83 & 5.66\\
				& & $\bn_2$ & 7.93 & 8.69 & 5.87 & 7.36 & 4.73 & 7.98 & 9.63 & 2.69 & 5.79 & 4.46 & 6.99 & 8.86 & 0.17 & 5.94 & 5.58\\
				& & $\bn_3$ & 6.49 & 7.90 & 4.79 & 5.58 & 4.72 & 7.20 & 9.92 & 1.92 & 5.13 & 4.21 & 6.64 & 8.86 & 0.38 & 4.92 & 4.87\\ \cline{2-18}
				& \multirow{3}{*}{500}& $\bn_1$& 8.90 & 11.04 & 7.65 & 11.24 & 5.14 & 7.74 & 9.80 & 4.21 & 15.25 & 5.62 & 7.50 & 9.18 & 0.01 & 6.80 & 4.81\\
				& & $\bn_2$ & 8.88 & 11.03 & 6.15 & 8.10 & 4.86 & 9.30 & 11.27 & 3.20 & 8.70 & 6.45 & 10.96 & 12.48 & 0.11 & 6.97 & 5.80\\
				& & $\bn_3$& 8.33 & 10.01 & 4.75 & 5.86 & 5.49 & 8.39 & 10.45 & 1.59 & 7.42 & 4.29 & 7.57 & 9.49 & 0.02 & 4.49 & 4.80\\ \hline
				\multirow{9}{*}{2} & \multirow{3}{*}{50} & $\bn_1$& 16.54 & 9.91 & 4.06 & 4.64 & 4.38 & 9.28 & 9.67 & 3.05 & 6.86 & 3.81 & 9.30 & 10.35 & 0.22 & 6.69 & 6.63\\
				&  & $\bn_2$ & 15.72 & 9.14 & 2.69 & 3.81 & 4.39 & 8.64 & 8.92 & 2.44 & 4.79 & 4.47 & 9.56 & 11.42 & 0.07 & 4.67 & 5.02\\
				&  & $\bn_3$& 14.44 & 8.74 & 2.95 & 3.39 & 4.28 & 10.21 & 10.75 & 2.34 & 4.64 & 4.68 & 8.84 & 9.87 & 0.18 & 4.51 & 4.67\\ \cline{2-18}
				& \multirow{3}{*}{100}& $\bn_1$& 12.49 & 9.66 & 2.82 & 4.51 & 4.38 & 10.61 & 11.86 & 2.85 & 8.33 & 6.55 & 8.94 & 10.24 & 0.02 & 6.36 & 5.80\\
				& & $\bn_2$ & 12.82 & 10.59 & 2.68 & 3.53 & 4.77 & 9.48 & 10.53 & 2.52 & 7.07 & 4.80 & 9.84 & 11.94 & 0.01 & 5.03 & 5.35\\
				& & $\bn_3$& 12.05 & 9.69 & 3.22 & 4.36 & 4.51 & 6.91 & 8.96 & 2.48 & 5.30 & 3.61 & 9.50 & 10.73 & 0.07 & 5.73 & 5.87\\ \cline{2-18}
				& \multirow{3}{*}{500}& $\bn_1$& 10.00 & 11.61 & 3.82 & 6.21 & 4.73 & 7.81 & 10.05 & 2.68 & 9.85 & 4.51 & 8.18 & 10.01 & 0.07 & 6.06 & 5.73\\
				& & $\bn_2$ & 9.01 & 10.23 & 3.25 & 4.82 & 4.48 & 8.79 & 10.39 & 1.99 & 6.31 & 4.63 & 7.73 & 9.18 & 0.08 & 4.72 & 4.46\\
				& & $\bn_3$& 9.39 & 10.24 & 2.86 & 3.64 & 5.25 & 9.39 & 11.03 & 2.23 & 7.06 & 6.18 & 8.57 & 9.91 & 0.07 & 4.73 & 5.67\\ \hline
				\multirow{9}{*}{3} & \multirow{3}{*}{50} & $\bn_1$& 37.22 & 9.75 & 1.20 & 1.85 & 5.38 & 12.26 & 11.28 & 2.46 & 5.93 & 5.40 & 8.42 & 9.95 & 0.28 & 4.02 & 5.33\\
				&  & $\bn_2$ & 43.28 & 11.73 & 1.26 & 1.76 & 5.51 & 12.09 & 11.18 & 0.73 & 3.45 & 4.36 & 9.58 & 11.61 & 0.27 & 5.62 & 5.73\\
				&  & $\bn_3$& 41.61 & 11.92 & 2.21 & 2.87 & 4.23 & 14.68 & 13.61 & 2.02 & 4.64 & 5.76 & 11.44 & 12.32 & 0.04 & 5.20 & 6.28\\  \cline{2-18}
				& \multirow{3}{*}{100}& $\bn_1$& 21.16 & 10.50 & 0.43 & 0.77 & 4.99 & 10.86 & 10.79 & 1.82 & 5.32 & 5.52 & 9.92 & 11.49 & 0.18 & 5.51 & 6.27\\
				& & $\bn_2$ & 19.50 & 11.65 & 0.79 & 1.15 & 4.78 & 9.92 & 11.22 & 1.18 & 4.78 & 4.69 & 9.97 & 11.11 & 0.24 & 4.14 & 5.78\\
				& & $\bn_3$& 22.79 & 13.01 & 1.22 & 1.57 & 4.74 & 10.61 & 11.28 & 1.58 & 5.04 & 5.61 & 10.32 & 12.06 & 0.02 & 4.46 & 5.35\\ \cline{2-18}
				& \multirow{3}{*}{500}& $\bn_1$& 8.86 & 9.01 & 0.48 & 1.12 & 4.33 & 8.19 & 9.52 & 1.08 & 5.28 & 4.72 & 7.34 & 8.57 & 0.00 & 4.71 & 4.60\\
				&& $\bn_2$ & 9.51 & 9.12 & 0.49 & 0.83 & 5.53 & 7.51 & 9.13 & 0.89 & 3.83 & 5.05 & 9.91 & 11.73 & 0.09 & 5.08 & 6.15\\
				& & $\bn_3$& 10.84 & 10.20 & 1.21 & 1.34 & 4.49 & 10.21 & 11.60 & 1.14 & 4.39 & 6.11 & 8.27 & 10.25 & 0.01 & 3.96 & 4.35\\ \hline
				\multicolumn{3}{c}{ARE}& 193.18 &99.04 &45.03& 45.27& 8.27& 84.78& 108.97& 52.90& 38.31& 14.11& 78.50& 110.61 &97.22&16.02& 13.16\\ \hline
		\end{tabularx}}
\end{table}

\begin{figure}
	\centering 
	\includegraphics[width=\textwidth]{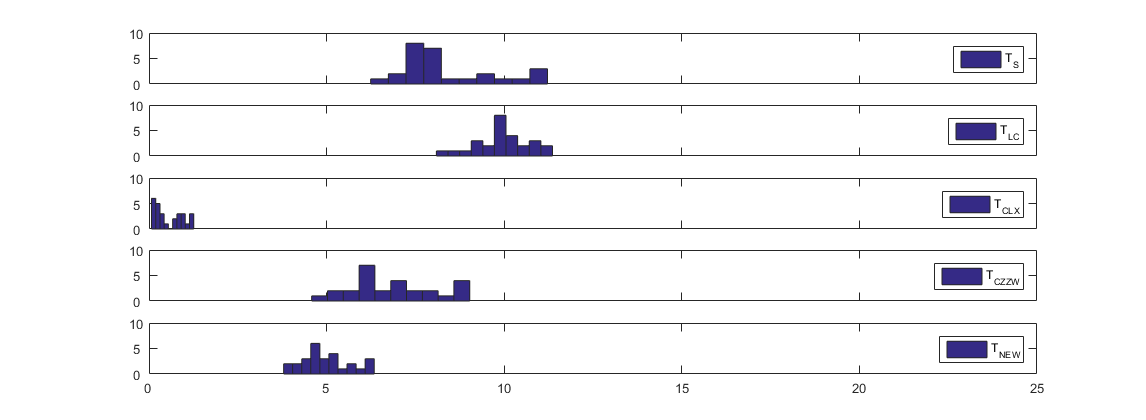}
	\caption{	\label{size1.fig} Histograms of the empirical sizes (in $\%$) of $T_{\TS}, T_{\LC}, T_{\CLX}, T_{\CZZW}$, and $T_{\NEW}$ (from top to bottom) in  Simulation 1. The nominal size is 5\%.}
\end{figure}

We first compare the empirical sizes of $T_{\TS}, T_{\LC}, T_{\CLX}, T_{\CZZW}$, and $T_{\NEW}$. Table~\ref{size1.tab} displays the empirical sizes of the five considered tests   with the last row displaying their ARE values associated the three values of $\rho$, from which we can draw several conclusions in terms of size control. First of all, $T_{\NEW}$ performs very well regardless of whether the data are  less correlated ($\rho=0.25$), moderately correlated ($\rho=0.5$), or highly correlated ($\rho=0.9$) since its empirical sizes range from $3.61\%$ to $6.63\%$ and its ARE values are $8.27$, $14.11$, and $13.16$ for $\rho=0.25,0.5$, and $0.9$, respectively.  Second, $T_{\TS}$ is very liberal with its empirical sizes ranging from $6.23\%$ to $43.28\%$ and its associated  ARE values being $193.18, 84.78, 78.50$, respectively. It is also seen that $T_{\TS}$ is much more liberal for skewed distributions (under Model 3) than for symmetric distributions (under Models 1 and 2). This is not a surprise since $T_{\TS}$ is proposed for normal data only. Third, $T_{\LC}$ is generally less liberal than $T_{\TS}$ since it is proposed for both normal and non-normal data. Nevertheless, it is still quite liberal with its empirical sizes ranging from $7.13\%$ to $13.01\%$ and its associated  ARE values being $99.04, 108.97, 110.61$,  respectively. Fourth,  $T_{\CLX}$ is very conservative with its empirical size as small as $0.00\%$,  especially when the data are not normally distributed or highly correlated and its associated ARE values being $45.03,52.90$, and $97.22$, respectively. Finally, $T_{\CZZW}$ does improve $T_{\CLX}$ with its associated  ARE values being $45.27,38.31$, and $16.02$, respectively. However it  can  also be very liberal with its empirical size being as large as $15.25\%$ and be very conservative with its empirical size being as small as $0.77\%$. Therefore,  Table~\ref{size1.tab} demonstrates that in terms of size control, the proposed normal-reference test $T_{\NEW}$  outperforms the four existing competitors $T_{\TS}, T_{\LC}, T_{\CLX}$, and $T_{\CZZW}$ substantially.  Some of the above conclusions may be further verified visually by Figure~\ref{size1.fig} which  displays  the histograms of the empirical sizes of $T_{\TS}, T_{\LC}, T_{\CLX}, T_{\CZZW}$, and $T_{\NEW}$ (from top to bottom).

To explain why $T_{\TS}$ and $T_{\LC}$ do not perform well in terms of size control, Table~\ref{df1.tab} displays the estimated approximate degrees of freedom $\hd$ (\ref{hbetadf.sec2})  of $T_{\NEW}$ under various settings in this simulation study.  It is seen that the values of $\hd$ all are rather  small. This  shows  that the null distribution of $T_{\NEW}$ is unlikely to be normal, and hence the normal approximations to the null distributions of $T_{\TS}$ and $T_{\LC}$ are unlikely to be adequate.  This partially explains why $T_{\TS}$ and $T_{\LC}$ are not accurate in size control regardless of whether the data are less correlated ($\rho=0.25$), moderately correlated ($\rho=0.5$), or highly correlated ($\rho=0.9$). Further,  it is seen that with the value of $\rho$ increasing, the value of  $\hd$ decreases. This means  that the more highly correlated the data  are, the less adequate the normal approximations to the null distributions of $T_{\TS}$ and $T_{\LC}$  would be.

\begin{table}[!h]
	\caption{\label{df1.tab}	 Estimated approximate degrees of freedom of $T_{\NEW}$ under various settings in Simulation 1.}
	\centering
		\setlength{\tabcolsep}{3pt}	
		\scalebox{0.9}{
			\begin{tabular}{ccccccccccc}
				\hline
 				\multicolumn{2}{c}{}        & \multicolumn{3}{c}{Model 1} & \multicolumn{3}{c}{Model 2} & \multicolumn{3}{c}{Model 3}  \\ \hline
				$p$ & $\bn$ & $\rho = 0.25$ & $\rho$ = 0.5  & $\rho$ = 0.9 & $\rho = 0.25$  & $\rho = 0.5$  & $\rho = 0.9$ & $\rho = 0.25$  & $\rho = 0.5$ & $\rho = 0.9$  \\ \hline
				\multirow{3}{*}{50}  & $\bn_1$   & 2.41  & 1.66  & 1.26  & 3.73  & 2.40   & 1.27 & 4.08  & 1.38  & 1.20   \\ \cline{2-11}
				&  $\bn_2$  & 2.93  & 1.19  & 1.06  & 3.96  & 1.17  & 1.08  &3.61  & 1.21  & 1.05  \\ \cline{2-11}
				&  $\bn_3$ & 2.52  & 1.09  & 1.02  & 1.90   & 1.08  & 1.02  &2.85  & 1.10   & 1.03  \\ \hline
				\multirow{3}{*}{100} & $\bn_1$   & 4.21  & 1.47  & 1.31  & 3.11  & 1.38  & 1.23 & 3.20  & 2.05  & 1.20   \\ \cline{2-11}
				& $\bn_2$  & 2.44  & 1.16  & 1.06  & 3.24  & 1.92  & 1.07 &2.65  & 1.16  & 1.06  \\ \cline{2-11}
				& $\bn_3$ & 1.45  & 1.06  & 1.02  & 2.51  & 1.07  & 1.02 & 1.94   & 1.08  & 1.02  \\ \hline
				\multirow{3}{*}{500} & $\bn_1$   & 6.85  & 1.40   & 1.23  & 5.97  & 1.72  & 1.37 & 1.95  & 1.51  & 1.24  \\ \cline{2-11}
				& $\bn_2$  & 1.86  & 1.09  & 1.05  & 3.23  & 1.09  & 1.05 & 1.61  & 1.09  & 1.05  \\ \cline{2-11}
				& $\bn_3$ & 1.36  & 1.04  & 1.03  & 1.45  & 1.04  & 1.03 & 1.29  & 1.05  & 1.02  \\ \hline
		\end{tabular}}
\end{table}

We now compare the empirical powers of  $T_{\TS}$, $T_{\LC}$, $T_{\CLX}$, $T_{\CZZW}$, and $T_{\NEW}$. Table~\ref{power1.tab} presents  the empirical powers of the five tests  under various configurations in Simulation 1. It is seen that  the empirical powers of the five tests when $\rho_1=\rho_2=0.5$ are generally smaller than their empirical powers when $\rho_1\neq \rho_2$. This is not a surprise since in this case,  the differences between $\bSigma_1-\bSigma_2$ just come from the differences between $\sigma_{1k}$'s and $\sigma_{2k}$'s.  Notice also  that the empirical sizes of the tests have strong impact on their empirical powers. Usually when  the empirical size of a test is larger than that of another test under a setting, its empirical power is often also larger than that of another test.  From Table~\ref{power1.tab},  we can see  that  the empirical powers of $T_{\TS}, T_{\LC}, $ and $T_{\CZZW}$ are generally ``larger'' than those of $T_{\NEW}$ and the empirical powers of $T_{\CLX}$ are generally ``smaller'' than those of $T_{\NEW}$. This is because  from Table~\ref{size1.tab} and Figure~\ref{size1.fig}, we can see that $T_{\TS}, T_{\LC}, $ and $T_{\CZZW}$ are generally  more liberal than $T_{\NEW}$,  $T_{\CLX}$ is more conservative than $T_{\NEW}$, while $T_{\NEW}$ has much better size control than the other tests.  Therefore, when a test does not have a good size control, its empirical powers can be misleading and hence it is very important to make a test have  a good size control.

\begin{table}
	\caption{\label{power1.tab}	 Empirical power (in $\%$)  of $T_{\TS}, T_{\LC}, T_{\CLX}, T_{\CZZW}$, and $T_{\NEW}$  in Simulation 1.}
	\centering
		\setlength{\tabcolsep}{3pt}	
		\scalebox{0.75}{
			\begin{tabularx}{1.3\textwidth}{@{\extracolsep{\fill}}cccccccccccccccccc }
				\hline
				\textbf{}& \textbf{}& \textbf{}& \multicolumn{5}{c}{$(\rho_1,\rho_2)= (0.5,0.5)$}& \multicolumn{5}{c}{$(\rho_1,\rho_2) = (0.5,0.25)$}& \multicolumn{5}{c}{$(\rho_1,\rho_2) = (0.5,0.9)$}\\ \hline
				Model & $p$ & $\bn$ & $T_{\TS}$ &  $T_{\LC}$ & $T_{\CLX}$ &  $T_{\CZZW}$ & $T_{\NEW}$ & $T_{\TS}$ &  $T_{\LC}$ & $T_{\CLX}$ &  $T_{\CZZW}$ & $T_{\NEW}$ & $T_{\TS}$ &  $T_{\LC}$ & $T_{\CLX}$ &  $T_{\CZZW}$ & $T_{\NEW}$ \\ \hline
				\multirow{9}[0]{*}{1} & \multirow{3}[0]{*}{50}    & $\bn_1$    & 8.30  & 9.54  & 4.37  & 9.34  & 4.31  & 80.69 & 82.17 & 34.40 & 50.61 & 63.88 & 71.51 & 74.46 & 61.83 & 85.47 & 66.95 \\
				&     &$\bn_2$   & 9.73  & 11.86 & 5.45  & 10.12 & 6.60  & 94.62 & 95.64 & 61.44 & 73.21 & 87.91 & 87.03 & 89.00 & 78.53 & 93.26 & 82.85 \\
				&     & $\bn_3$   & 12.59 & 15.81 & 7.76  & 15.31 & 7.24  & 100.00 & 100.00 & 99.14 & 99.54 & 100.00 & 99.77 & 99.74 & 100.00 & 100.00 & 99.62 \\
				\cline{2-18}
				& \multirow{3}[0]{*}{100}    & $\bn_1$    & 8.69  & 10.20 & 4.95  & 11.51 & 5.45  & 82.53 & 85.38 & 38.21 & 55.87 & 66.69 & 72.57 & 75.84 & 66.16 & 89.64 & 68.63 \\
				&   & $\bn_2$    & 9.22  & 11.61 & 4.47  & 10.58 & 5.42  & 93.56 & 94.46 & 60.70 & 75.16 & 87.30 & 87.92 & 89.58 & 80.41 & 95.94 & 84.84 \\
				&   & $\bn_3$    & 11.82 & 13.28 & 6.78  & 16.96 & 6.48  & 99.96 & 99.97 & 98.73 & 99.55 & 99.93 & 100.00 & 100.00 & 99.94 & 100.00 & 99.98 \\
				\cline{2-18}
				& \multirow{3}[0]{*}{500}    & $\bn_1$     & 9.92  & 11.92 & 5.67  & 16.66 & 7.02  & 80.79 & 82.73 & 37.09 & 60.60 & 64.27 & 71.13 & 74.88 & 62.24 & 93.10 & 68.28 \\
				&   & $\bn_2$    & 9.75  & 12.30 & 4.38  & 14.88 & 6.09  & 95.26 & 96.57 & 62.83 & 82.32 & 88.93 & 91.03 & 91.94 & 81.67 & 97.73 & 87.00 \\
				&   & $\bn_3$    & 11.99 & 14.13 & 7.17  & 21.34 & 6.79  & 100.00 & 100.00 & 99.80 & 100.00 & 100.00 & 99.97 & 99.79 & 100.00 & 100.00 & 99.90 \\
				\hline
				\multirow{9}[0]{*}{2} & \multirow{3}[0]{*}{50}     & $\bn_1$    & 9.45  & 10.43 & 3.68  & 7.79  & 5.64  & 81.02 & 81.27 & 33.00 & 47.94 & 59.38 & 71.06 & 72.93 & 55.00 & 83.73 & 63.02 \\
				&    & $\bn_2$    & 12.30 & 13.08 & 3.84  & 9.14  & 6.83  & 94.38 & 94.64 & 54.83 & 68.94 & 83.95 & 87.78 & 88.33 & 77.52 & 92.21 & 82.50 \\
				&   & $\bn_3$    & 15.78 & 16.83 & 6.75  & 14.19 & 7.56  & 100.00 & 99.92 & 98.24 & 99.18 & 99.92 & 99.70 & 99.82 & 99.30 & 99.90 & 99.21 \\
				\cline{2-18}
				& \multirow{3}[0]{*}{100}    & $\bn_1$     & 10.34 & 11.24 & 3.30  & 9.50  & 5.06  & 82.38 & 83.13 & 35.49 & 51.66 & 65.47 & 68.50 & 72.37 & 55.66 & 85.80 & 64.26 \\
				&   & $\bn_2$    & 9.91  & 11.22 & 4.02  & 9.01  & 5.54  & 95.96 & 96.07 & 61.08 & 75.18 & 87.43 & 87.59 & 89.40 & 75.17 & 94.24 & 82.74 \\
				&    & $\bn_3$    & 12.06 & 14.83 & 6.00  & 14.68 & 7.19  & 99.95 & 99.92 & 98.67 & 99.43 & 99.89 & 99.84 & 99.81 & 99.45 & 99.93 & 99.55 \\
				\cline{2-18}
				& \multirow{3}[0]{*}{500}   & $\bn_1$     & 8.79  & 10.84 & 3.92  & 12.89 & 6.06  & 85.27 & 86.48 & 31.45 & 57.33 & 68.16 & 70.60 & 74.25 & 53.61 & 89.15 & 66.48 \\
				&   & $\bn_2$  & 11.09 & 12.45 & 2.86  & 9.94  & 6.07  & 96.12 & 96.87 & 63.04 & 81.11 & 90.84 & 88.73 & 90.11 & 75.46 & 95.13 & 85.33 \\
				&  & $\bn_3$    & 12.85 & 15.94 & 5.81  & 20.98 & 6.99  & 100.00 & 100.00 & 98.67 & 99.33 & 99.92 & 100.00 & 100.00 & 100.00 & 100.00 & 100.00 \\
				\hline
				\multirow{9}[0]{*}{3} & \multirow{3}[0]{*}{50}     & $\bn_1$     & 11.58 & 10.26 & 2.76  & 6.62  & 4.87  & 82.01 & 79.00 & 35.85 & 52.51 & 53.23 & 72.93 & 73.98 & 51.55 & 79.24 & 63.77 \\
				&    & $\bn_2$    & 14.83 & 14.05 & 3.15  & 7.91  & 6.06  & 94.94 & 93.32 & 63.50 & 77.46 & 81.56 & 89.10 & 89.56 & 73.63 & 93.13 & 81.01 \\
				&   & $\bn_3$    & 16.08 & 15.52 & 4.65  & 12.90 & 7.07  & 100.00 & 100.00 & 98.43 & 99.40 & 99.71 & 99.83 & 99.81 & 99.54 & 99.94 & 99.35 \\
				\cline{2-18}
				& \multirow{3}[0]{*}{100}   & $\bn_1$    & 12.21 & 12.13 & 3.17  & 9.23  & 7.21  & 82.56 & 82.29 & 32.65 & 51.66 & 59.68 & 73.08 & 76.39 & 53.31 & 83.43 & 69.29 \\
				&   & $\bn_2$    & 13.25 & 12.75 & 2.66  & 7.13  & 5.94  & 95.35 & 94.93 & 61.27 & 77.40 & 84.86 & 89.02 & 89.80 & 72.61 & 93.25 & 83.22 \\
				&    & $\bn_3$   & 13.84 & 14.12 & 4.96  & 10.53 & 6.07  & 99.93 & 99.97 & 99.39 & 99.50 & 99.67 & 99.88 & 99.80 & 99.67 & 99.97 & 99.50 \\
				\cline{2-18}
				& \multirow{3}[0]{*}{500}   & $\bn_1$     & 9.83  & 11.58 & 1.86  & 7.54  & 6.48  & 82.48 & 83.35 & 26.46 & 52.44 & 65.87 & 71.09 & 74.06 & 43.30 & 85.47 & 67.67 \\
				&   & $\bn_2$    & 11.06 & 13.34 & 2.02  & 9.15  & 6.31  & 95.19 & 95.00 & 55.77 & 78.84 & 88.55 & 88.29 & 89.64 & 64.20 & 94.13 & 84.63 \\
				&   & $\bn_3$    & 13.14 & 15.07 & 4.47  & 16.33 & 7.33  & 100.00 & 100.00 & 99.28 & 99.95 & 100.00 & 100.00 & 100.00 & 99.89 & 100.00 & 99.91 \\
				\hline
		\end{tabularx}}
\end{table}

\subsection{Simulation 2}\label{sim2.subsec}

In this simulation study, we continue to compare $T_{\NEW}$ against  $T_{\TS}$, $T_{\LC}$, $T_{\CLX}$, and $T_{\CZZW}$ for the two-sample equal-covariance matrix testing  problem (\ref{covtest.sec2}) but with the two samples (\ref{twosamp.sec1}) generated from the following moving average model:
	\[
	y_{ijk}=z_{ijk}+\theta_{i1}z_{ij(k+1)}+\cdots+\theta_{im_{i}}z_{ij(k+m_{i})},k=1,...,p;\;j=1,...,n_i;\;i=1,2,
	\]
	where $y_{ijk}$ denotes the $k$th component of $\by_{ij}$, and $z_{ij\ell},\ell=1,\ldots,p+m_i;\;i=1,2$ are i.i.d. random variables generated in the same way as described in Simulation 1. The covariance matrix difference $\bSigma_1-\bSigma_2$ is then  determined by  $m_i,i=1,2$ and $\theta_{ij},j=1,\ldots,m_i;\;i=1,2$. When $m_1=m_2=m$ and $\theta_{1j}=\theta_{2j}=\theta_j,j=1,\ldots,m$, we have $\bSigma_1=\bSigma_2$ so that the generated two samples (\ref{twosamp.sec1}) satisfy  the null hypothesis in (\ref{covtest.sec2}). For size comparison, we set $m=0.5p$, and let $\theta_j, j=1,\ldots,m$ be  generated from $U[2,3]$, the uniform distribution over $[2,3]$. For power comparison, we set $m_1=0.5p$, $m_2=0.4p$  and let  $\theta_{1j},j=1,\ldots,m_1$ be generated from $U[2,3]$, and $\theta_{2j},j=1,\ldots,m_2$ be  generated from  $U[3,4]$. The tuning parameters $p$ and $\bn$ are specified in the same way as in Simulation 1.

\begin{figure}[!h]
	\centering 
	\includegraphics[width=\textwidth]{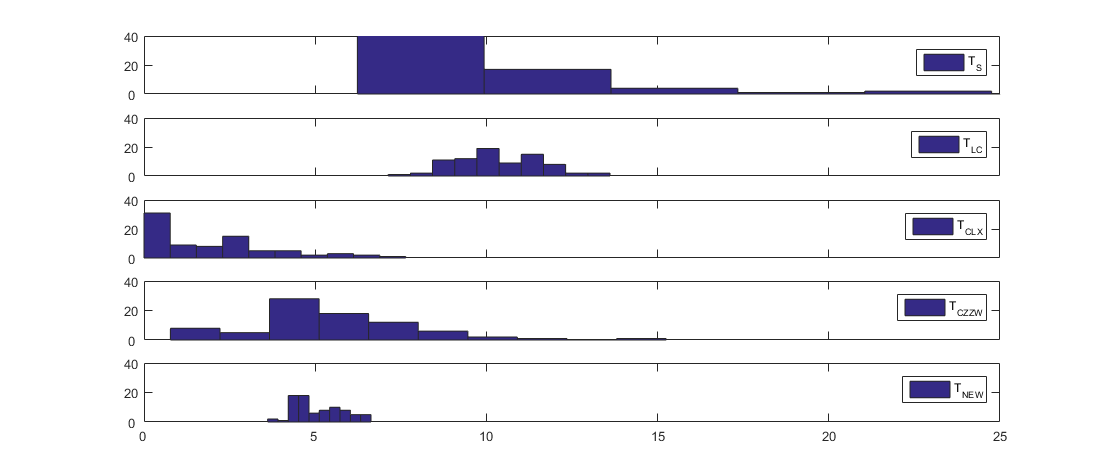}
	\caption{	\label{size2.fig}Histograms of the empirical sizes (in $\%$) of $T_{\TS}, T_{\LC}, T_{\CLX}, T_{\CZZW}$, and $T_{\NEW}$ (from top to bottom) in  Simulation 2.}

\end{figure}

Table~\ref{size2.tab} displays the empirical sizes and powers of $T_{\TS}$, $T_{\LC}$, $T_{\CLX}$, $T_{\CZZW}$, and $T_{\NEW}$ in this simulation study.  In terms of size control, it is seen that the conclusions drawn from this table and Figure~\ref{size2.fig} are similar to those drawn from Table~\ref{size1.tab} and Figure~\ref{size1.fig}. That is,  $T_{\NEW}$ performs well with its empirical sizes ranging from $3.79\%$ to $6.34\%$. It performs much better than $T_{\TS}, T_{\LC}$, and $T_{\CZZW}$ which  are generally quite liberal with their empirical sizes ranging from $6.24\%$ to $11.22\%$, $8.09\%$ to $11.36\%$, and $4.58\%$ to $9.03\%$, respectively. It also performs much better than $T_{\CLX}$ which is generally quite conservative with their empirical sizes ranging from $0.07\%$ to $1.26\%$. These conclusions are also indicated by  the ARE values of  $T_{\TS}$, $T_{\LC}$, $T_{\CLX}$, $T_{\CZZW}$, and $T_{\NEW}$ as listed in the last row of the table.  In terms of powers,
it is seen that the conclusions drawn from  this table are similar to those drawn from Table~\ref{power1.tab}. That is, the empirical powers of $T_{\TS}, T_{\LC}, $ and $T_{\CZZW}$ are generally ``larger'' than those of $T_{\NEW}$ while the empirical powers of $T_{\CLX}$ are generally ``smaller'' than those of $T_{\NEW}$. This is also  because  from Table~\ref{size2.tab} and Figure~\ref{size2.fig}, we can see that $T_{\TS}, T_{\LC}, $ and $T_{\CZZW}$ are generally  more liberal than $T_{\NEW}$,  $T_{\CLX}$ is more conservative
than $T_{\NEW}$, while $T_{\NEW}$ has much better size control than the other tests.

\begin{table}
	\caption{\label{size2.tab}	Empirical sizes and powers (in $\%$) of $T_{\TS}$, $T_{\LC}$, $T_{\CLX}$, $T_{\CZZW}$, and $T_{\NEW}$ in Simulation 2.}
	\centering
		\setlength{\tabcolsep}{3pt}	
		\scalebox{0.9}{
			\begin{tabularx}{\textwidth}{@{\extracolsep{\fill}}ccccccccccccc}
				\hline
				\multicolumn{3}{c}{} &\multicolumn{5}{c}{$\bSigma_1=\bSigma_2$} &\multicolumn{5}{c}{$\bSigma_1\neq\bSigma_2$} \\
				\hline
				Model & $p$ & $\bn$ &$T_{\TS}$ & $T_{\LC}$ & $T_{\CLX}$ & $T_{\CZZW}$ & $T_{\NEW}$ &$T_{\TS}$ & $T_{\LC}$ & $T_{\CLX}$ & $T_{\CZZW}$ & $T_{\NEW}$ \\ \hline
				\multirow{9}{*}{1} & \multirow{3}{*}{50} & $\bn_1$ &6.76 & 9.86 & 1.26 & 8.60 & 5.55 & 38.93 & 46.50 & 15.34 & 45.35 & 38.12\\
				& & $\bn_2$ &7.76 & 9.44 & 1.14 & 7.75 & 6.14 & 59.79 & 65.95 & 28.78 & 66.16 & 57.30\\
				& & $\bn_3$ &7.60 & 10.14 & 1.01 & 5.47 & 4.87 & 99.01 & 99.55 & 88.36 & 99.17 & 98.10\\
				\cline{2-13}
				&  \multirow{3}{*}{100} & $\bn_1$ &7.57 & 9.78 & 0.84 & 7.63 & 5.65 & 46.33 & 53.59 & 13.86 & 52.67 & 45.63\\
				& & $\bn_2$ &7.49 & 9.77 & 0.91 & 7.92 & 5.28 & 70.93 & 77.50 & 29.25 & 74.45 & 67.84\\
				& & $\bn_3$ &7.25 & 9.47 & 0.30 & 5.24 & 3.97 & 99.50 & 99.80 & 91.19 & 99.60 & 99.16\\
				\cline{2-13}
				& \multirow{3}{*}{500} & $\bn_1$ &8.20 & 10.74 & 0.28 & 7.04 & 6.10 & 52.62 & 60.63 & 5.00 & 54.82 & 51.83\\
				& & $\bn_2$ &7.29 & 9.21 & 0.07 & 6.58 & 5.03 & 77.39 & 82.67 & 11.38 & 75.73 & 75.27\\
				& & $\bn_3$ &7.41 & 9.75 & 0.16 & 6.50 & 4.66 & 100.00 & 100.00 & 76.76 & 99.96 & 100.00\\
				\hline
				\multirow{9}{*}{2} & \multirow{3}{*}{50} & $\bn_1$ &10.64 & 10.53 & 0.84 & 8.61 & 6.34 & 40.50 & 44.47 & 12.95 & 43.40 & 32.92\\
				& & $\bn_2$ &9.74 & 11.13 & 0.76 & 6.16 & 4.64 & 64.13 & 67.36 & 25.37 & 60.61 & 54.18\\
				& & $\bn_3$ &9.44 & 10.68 & 0.95 & 6.11 & 5.26 & 98.33 & 98.58 & 82.20 & 97.41 & 96.40\\
				\cline{2-13}
				&  \multirow{3}{*}{100} & $\bn_1$ &7.72 & 9.73 & 0.79 & 8.29 & 4.91 & 48.43 & 53.99 & 13.06 & 50.41 & 44.22\\
				& & $\bn_2$ &8.05 & 9.36 & 0.29 & 5.97 & 4.30 & 71.25 & 75.87 & 26.28 & 70.24 & 65.42\\
				& & $\bn_3$ &7.84 & 9.96 & 0.35 & 7.19 & 4.67 & 99.71 & 99.83 & 87.26 & 99.19 & 99.20\\
				\cline{2-13}
				& \multirow{3}{*}{500} & $\bn_1$ &7.85 & 9.01 & 0.25 & 9.03 & 5.14 & 50.23 & 56.16 & 3.54 & 53.56 & 47.83\\
				& & $\bn_2$ &6.24 & 8.09 & 0.10 & 6.19 & 3.79 & 77.64 & 82.90 & 11.27 & 76.49 & 74.47\\
				& & $\bn_3$ &7.93 & 10.23 & 0.12 & 6.15 & 5.64 & 100.00 & 100.00 & 77.90 & 100.00 & 100.00\\
				\hline
				\multirow{9}{*}{3} & \multirow{3}{*}{50} & $\bn_1$ &11.22 & 10.92 & 1.16 & 6.99 & 4.11 & 44.28 & 44.41 & 11.51 & 39.56 & 31.00\\
				& & $\bn_2$ &11.19 & 11.36 & 1.17 & 7.52 & 4.41 & 63.34 & 63.07 & 22.38 & 53.08 & 48.00\\
				& & $\bn_3$ &11.09 & 10.30 & 0.73 & 4.58 & 4.65 & 98.35 & 98.24 & 75.27 & 95.16 & 93.55\\
				\cline{2-13}
				&  \multirow{3}{*}{100} & $\bn_1$ &7.26 & 8.48 & 0.33 & 7.13 & 4.41 & 51.38 & 54.85 & 10.60 & 50.17 & 44.30\\
				& & $\bn_2$ &8.77 & 9.77 & 0.47 & 5.80 & 4.79 & 73.18 & 76.68 & 21.58 & 68.52 & 63.70\\
				& & $\bn_3$ &9.38 & 11.02 & 0.33 & 5.71 & 4.76 & 99.40 & 99.65 & 84.22 & 98.60 & 98.15\\
				\cline{2-13}
				& \multirow{3}{*}{500} & $\bn_1$ &8.50 & 10.13 & 0.16 & 8.64 & 6.03 & 51.62 & 57.00 & 5.14 & 54.11 & 48.73\\
				& & $\bn_2$ &7.04 & 9.26 & 0.28 & 6.01 & 4.34 & 77.63 & 81.95 & 12.46 & 76.24 & 74.44\\
				& & $\bn_3$ &7.96 & 9.92 & 0.13 & 6.22 & 5.13 & 99.70 & 99.86 & 76.31 & 99.52 & 99.50\\
				\hline
				\multicolumn{3}{c}{ARE}&66.81& 98.55& 88.76& 37.68& 11.12\\
				\hline
		\end{tabularx}}
\end{table}

\section{Application to the Financial Data}\label{appl.sec}


In this section, we apply $T_{\NEW}$, together with the other four tests, i.e., $T_{\TS}$, $T_{\LC}$, $T_{\CLX}$, and $T_{\CZZW}$ to the financial data set briefly described in Section~\ref{intro.sec}. It contains the daily PD  of  $235$ companies from the mainland of China (PD-CHN) and $153$ companies from Hong Kong (PD-HKG) in 2007--2008. We  aim to check whether PD-CHN and PD-HKG have the same covariance matrix.

\begin{table}[!h]
	\caption{\label{vir.tab1} Testing the equality of covariance matrices of the financial data}
	\centering
		\begin{tabular}{llccc}
			\hline
			Hypothesis &Method & Statistic& $p$-value & d.f. \\ \hline
			\multirow{5}{*}{PD-CHN vs. PD-HKG}& $T_{\TS}$ & 43.33 & 0& \\ 
			\cline{2-5}
			& $T_{\LC}$ & 35.00 & 0& \\ 
			\cline{2-5}
			& $T_{\CLX}$& 34.21& $8.12\times 10^{-4}$ & \\ \cline{2-5}
			& $T_{\CZZW}$ & 5.85& 0& \\ \cline{2-5}
			& $T_{\NEW}$& 11.11 & $3.65 \times 10^{-5}$ & 1.07 \\ \hline
		\end{tabular}
\end{table}

Table~\ref{vir.tab1} presents  the testing results of applying the five considered tests to the financial data set. It is seen that the $p$-values of the five tests are essentially 0, suggesting that the covariance matrices for these two groups are significantly different. In addition, the estimated approximate degrees of freedom (d.f.) of $T_{\NEW}$ is around $1$, showing that the normal approximations to the null distributions  of $T_{\TS}$ and $T_{\LC}$ are unlikely to be adequate.

To further demonstrate the finite-sample performance  of $T_{\NEW}$ in terms of accuracy, we use the data set to calculate the empirical sizes of the five tests under consideration  for PD-CHN and PD-HKG, respectively. The empirical sizes are obtained from $10,000$ runs. In each run, for PD-CHN, we randomly split the $235$ companies into two subgroups, namely PD-CHN-1 and PD-CHN-2, with their sample sizes being 117 and 118, respectively. We then compute the $p$-values of the five tests to check the equality of covariance matrices between PD-CHN-1 and PD-CHN-2. The empirical size is calculated as the proportion of times that the $p$-values are smaller than the nominal level $\alpha=5\%$ based on the $10,000$ simulation runs. Similarly, for PD-HKG, in each run, we randomly split the $153$ companies into two subgroups, namely PD-HKG-1 and PD-HKG-2, with their sample sizes being 76 and 77, respectively. Then we compute the $p$-values of the five tests  to check the equality of covariance matrices between PD-HKG-1 and PD-HKG-2, and compute the empirical sizes of  the five considered tests. The empirical sizes of $T_{\TS}$, $T_{\LC}$, $T_{\CLX}$, $T_{\CZZW}$, and $T_{\NEW}$ are presented in Table~\ref{vir.tab2}. From this table, we can see that $T_{\NEW}$ has a much better level accuracy than the other four tests. Conclusions from the above two small-scale simulation studies based on  the financial data are consistent with those drawn from the two simulation studies presented  in Section~\ref{simul.sec}.
\begin{table}
	\caption{\label{vir.tab2} Empirical sizes ($\%$) of the financial data}
	\centering
		\begin{tabular}{lccccc}
			\hline
			Population &$T_{\TS}$ &$T_{\LC}$ &$T_{\CLX}$ &$T_{\CZZW}$ &$T_{\NEW}$ \\ \hline
			PD-CHN & 24.79 & 13.89 & 0.00 & 7.72 & 6.10\\ \cline{2-6}
			PD-HKG & 15.54 & 14.57 & 0.00 & 5.71 & 5.48 \\ \hline
		\end{tabular}
\end{table}

\section{Concluding Remarks}\label{remark.sec}

It is often of interest to test whether two high-dimensional samples have the same covariance matrix. Existing tests are either too liberal or too conservative in terms of size control. In this paper,  we propose and study a normal-reference test which has much better size control than several existing competitors as demonstrated by two simulation studies and a real data example.

In theoretical development, we assume that the mean vectors are zero to construct the test statistic while  in the implementation of the proposed test,  we actually   replace the original two samples with  the two samples obtained after subtracting the sample mean vectors. In the settings of  our two simulation studies, it seems that the effect of this replacement  is ignorable. However, it is still unclear if this replacement may strongly affect the performance of the proposed test in other settings. Further study in this direction is interesting and warranted.

In the implementation of the proposed test, we estimate the key quantities $\calK_2(T_{n,p,0}^*)$ and $\calK_3(T_{n,p,0}^*)$  via estimating $\tr(\bOmega_1^2),\tr(\bOmega_1\bOmega_2),\tr(\bOmega_2^2),\tr(\bOmega_1^3),\tr(\bOmega_2^3),\tr(\bOmega_1^2\bOmega_2)$, and $\tr(\bOmega_1\bOmega_2^2)$, using the two induced samples (\ref{newsamp.sec2}). Whereas, it is also possible to estimate $\tr(\bOmega_1^2),\tr(\bOmega_2^2)$, $\tr(\bOmega_1\bOmega_2),\tr(\bOmega_1^3),\tr(\bOmega_1^2\bOmega_2),\tr(\bOmega_1\bOmega_2^2)$, and $\tr(\bOmega_2^3)$ using the original two samples (\ref{twosamp.sec1}). It is of interest to know which approach is better. Further study in this direction is interesting and warranted.

\section*{Funding and Acknowledgement}

The work was  financially  supported by the National University of Singapore academic research grant R-155-000-121-114.

\appendix
\bigskip{}
\centerline{APPENDIX: Technical Proofs} \setcounter{equation}{0}
\global\long\def\theequation{A.\arabic{equation}}

\noindent Proof of Theorem~\ref{Tnexp.thm}. Notice that we can write $\tT_{n,p,0}$ and $\tT_{n,p,0}^*$ as  generalized quadratic forms as defined in \citet[Section S.2 of the Appendix]{wangxu2021}:
\[
	\tT_{n,p,0}=\sum_{1\le k<l\le n} a_{kl} \bxi_k^\top\bxi_l,\quad\mbox{ and }\quad
	\tT_{n,p,0}^*=\sum_{1\le k<l\le n} a_{kl} \bxi_k^{*\top}\bxi_l^*,
	\]
	where $\bxi_k=\bu_{1k}, \bxi_k^*=\bu_{1k}^*, k=1,\ldots,n_1$, $\bxi_{n_1+k}=\bu_{2k}, \bxi_{n_1+k}^*=\bu_{2k}^*, k=1,\ldots,n_2$,  and
	\[
	a_{kl}=\left\{
	\begin{array}{ll} \frac{2}{n_1(n_1-1)\sigma_T}, & \mbox{when $1\le k<l\le n_1$},\\
		\frac{2}{n_2(n_2-1)\sigma_T}, & \mbox{when $n_1+1\le k<l\le n$},\\
		-\frac{2}{n_1n_2\sigma_T},    & \mbox{when $1\le k\le n_1, n_1+1\le l\le n$}.
	\end{array}
	\right.
	\]
	We will employ Theorem S.1  of  \cite{wangxu2021} in the following proofs. To employ Theorem S.1  of  \cite{wangxu2021}, we need to check Assumptions S.1 and S.2 of  \cite{wangxu2021} first. Note that we have $\sigma_{kl}^2=\E\left(a_{kl}\bxi_k^\top\bxi_l\right)^2=a_{kl}^2\tr(\bOmega_k\bOmega_l), 1\le k<l\le n$ where $\bOmega_k=\bOmega_1$ when $1\le k\le n_1$, and $\bOmega_k=\bOmega_2$ when $n_1+1\le k\le n$.
	Under Condition C2, we have
	\[
	\begin{split}
		\E\left(a_{kl}\bxi_k^\top\bxi_l\right)^4&=a_{kl}^4\E(\bxi_k^\top\bxi_l)^4=a_{kl}^4\E\left[\E\{(\bxi_k^\top\bxi_l)^4|\bxi_k\}\right]\\
		&\le a_{kl}^4\E\left[\gamma\E^2\{(\bxi_k^\top\bxi_l)^2|\bxi_k\}\right]=\gamma a_{kl}^4\E\left\{(\bxi_k^\top\bOmega_l\bxi_k)^2\right\}\\
		&=\gamma a_{kl}^4\E\|\bOmega_l^{1/2}\bxi_k\|^4
		\le \gamma^2  a_{kl}^4\E^2(\|\bOmega_l^{1/2}\bxi_k\|^2)\\
		&=\gamma^2 a_{kl}^4\tr^2(\bOmega_k\bOmega_l)=\gamma^2 \sigma_{kl}^4.
	\end{split}
	\]
	Similarly, we can show that
	\[
	\E\left(a_{kl}\bxi_k^\top\bxi_l^*\right)^4\le\gamma^2 \sigma_{kl}^4,\quad
	\E\left(a_{kl}\bxi_k^{*\top}\bxi_l\right)^4\le\gamma^2 \sigma_{kl}^4,\quad \mbox{ and }\quad
	\E\left(a_{kl}\bxi_k^{*\top}\bxi_l^*\right)^4\le\gamma^2 \sigma_{kl}^4.
	\]
	In addition, the other conditions in  Assumptions S.1 and S.2 of \citet{wangxu2021} are also satisfied by $\tT_{n,p,0}$ and $\tT_{n,p,0}^*$, and  they are independent from each other.   Applying Theorem S.1 of \cite{wangxu2021}, we have
	\[
	\left\|\calL(\tT_{n,p,0})-\calL(\tT_{n,p,0}^*)\right\|_3\le \frac{\gamma^{3/2}}{3^{1/4}}\sum_{k=1}^n \sinf_k^{3/2}.
	\]
	By the Cauchy--Schwarz inequality, we have
	\[
	\sum_{k=1}^n \sinf_k^{3/2}\le \left(\sum_{k=1}^n \sinf_k\right)^{1/2}\left(\sum_{k=1}^n \sinf_k^2\right)^{1/2}=\left\{\left(\sum_{k=1}^n \sinf_k\right)\left(\sum_{k=1}^n \sinf_k^2\right)\right\}^{1/2}.
	\]
	Now, when $1\le k\le n_1$, by \citet[p.23 of  the Appendix]{wangxu2021}, we have
	\[
	\begin{split}
		\sinf_k&=\sum_{l=1}^{k-1}\sigma_{lk}^2+\sum_{l=k+1}^{n}\sigma_{kl}^2=\sum_{l=1}^{k-1}\sigma_{lk}^2+\sum_{l=k+1}^{n_1}\sigma_{kl}^2+\sum_{l=n_1+1}^{n}\sigma_{kl}^2\\
		&=\sum_{l=1}^{k-1}a_{lk}^2\tr(\bOmega_1^2)+\sum_{l=k+1}^{n_1}a_{kl}^2\tr(\bOmega_1^2)+\sum_{l=n_1+1}^{n}a_{kl}^2\tr(\bOmega_1\bOmega_2)\\
		&=\frac{4\tr(\bOmega_1^2)}{n_1^2(n_1-1)\sigma_T^2}+\frac{4\tr(\bOmega_1\bOmega_2)}{n_1^2n_2\sigma_T^2}
		=\frac{2G_1}{n_1},
	\end{split}
	\]
	where $G_1=2\{n_1(n_1-1)\}^{-1}\tr(\bOmega_1^2)/\sigma_T^2+2(n_1n_2)^{-1}\tr(\bOmega_1\bOmega_2)/\sigma_T^2$.
	Similarly, when $n_1+1\le k \le n$, we have
	$\sinf_k=2G_2/n_2$ where $G_2=2(n_1n_2)^{-1}\tr(\bOmega_1\bOmega_2)/\sigma_T^2+2\{n_2(n_2-1)\}^{-1}\tr(\bOmega_2^2)/\sigma_T^2$.  It is easy to see that  $0<G_1,G_2<1$ and $G_1+G_2=1$.
	Therefore, we have
	\[
	\begin{split}
		\sum_{k=1}^n\sinf_k&=\sum_{k=1}^{n_1} \sinf_k+\sum_{k=n_1+1}^n \sinf_k=2G_1+2G_2=2,\;\mbox{ and }\\
		\sum_{k=1}^n\sinf_k^2&=\sum_{k=1}^{n_1} \sinf_k^2+\sum_{k=n_1+1}^n  \sinf_k^2=\frac{4G_1^2}{n_1}+\frac{4G_2^2}{n_2}\le 4\left(\frac{1}{n_1}+\frac{1}{n_2}\right).
	\end{split}
	\]
	It follows that
	\[
	\sum_{k=1}^n \sinf_k^{3/2}\le 2^{3/2}\left(\frac{1}{n_1}+\frac{1}{n_2}\right)^{1/2}.
	\]
	Thus, we have
	\[
	\left\|\calL(\tT_{n,p,0})-\calL(\tT_{n,p,0}^*)\right\|_3\le \frac{(2\gamma)^{3/2}}{3^{1/4}} \left(\frac{1}{n_1}+\frac{1}{n_2}\right)^{1/2}.
	\]
	The proof is complete.

\noindent Proof of Theorem~\ref{Tnlim.thm}. Since for $i=1,2$, we have $\bu_{ij}^*, j=1,\ldots, n_i\iidsim N_{p^2}(\bzero,\bOmega_i)$, then $\barbu_1^*-\barbu_2^*\sim N_{p^2}(\bzero,\bOmega_n)$, and we can express $\|\barbu_1^*-\barbu_2^*\|^2=\bxi_{p^2}^{\top}\bOmega_n\bxi_{p^2}$, where $\bxi_{p^2}\sim N_{p^2}(\bzero,\bI_{p^2})$. Recall that
	\[
	T_{n,p,0}^*=\|\barbu_1^*-\barbu_2^*\|^2-\tr(\hbOmega_n^*),
	\]
	where $\hbOmega_i^*=(n_i-1)^{-1}\sum_{j=1}^{n_i}(\bu_{ij}^*-\barbu_i^*)(\bu_{ij}^*-\barbu_i^*)^\T,i=1,2$, and $\hbOmega_n^*=\hbOmega_1^*/n_1+\hbOmega_2^*/n_2$. It follows that
	\[
	\tr(\hbOmega_n^*)=\tr(\hbOmega_1^*)/n_1+\tr(\hbOmega_2^*)/n_2.
	\]
	Further, we have $\E\{\tr(\hbOmega_n^*)\}=\tr(\bOmega_n)$ and
	$\Var\{\tr(\hbOmega_n^*)\}=2\{n_1^2(n_1-1)\}^{-1}\tr(\bOmega_1^2)+2\{n_2^2(n_2-1)\}^{-1}\tr(\bOmega_2^2)$. Under Condition C1, as $n\to\infty$, we have
	\[
	\Var\{\tr(\hbOmega_n^*)/\tr(\bOmega_n)\}=\frac{2\tr(\bOmega_1^2)/\{n_1^2(n_1-1)\}+2\tr(\bOmega_2^2)/\{n_2^2(n_2-1)\}}{\tr^2(\bOmega_n)}\to 0,
	\]
	uniformly for all $p$. Thus, $\tr(\hbOmega_n^*)/\tr(\bOmega_n)\to 1$ in probability uniformly for all $p$. By (\ref{cumulants.sec2}), we have $\sigma_T^2=2\tr(\bOmega_n^2)\{1+\lito(1)\}$.  It follows that we have
	\[
	\tT_{n,p,0}^*=\frac{T_{n,p,0}^*}{\sigma_T}=\frac{\|\barbu_1^*-\barbu_2^*\|^2-\tr(\hbOmega_n^*)}{\sqrt{\sigma_T^2}}\\
	=\frac{\bxi_{p^2}^{\top}\bOmega_n\bxi_{p^2}-\tr(\bOmega_n)}{\sqrt{2\tr(\bOmega_n^2)}}\{1+\lito_p(1)\}.
	\]
	Under Condition C3, the expression (\ref{Tnlim.sec2})   follows from Corollary 1 of \cite{wangxu2021} immediately. The proof is complete.

\noindent Proof of Theorem~\ref{Tnpower.thm}.
	By (\ref{Tdecomp.sec2}) and under the local alternative  (\ref{powercond1.sec2}), we have \[ T_{n,p}=\left[T_{n,p,0}+\tr\left\{(\bSigma_1-\bSigma_2)^2 \right\}\right]\{1+\lito_p(1)\}.\]
	By (\ref{Tnsig2.sec2}) and (\ref{Omega.sec2}),  we have $\sigma_T^2=2\tr(\bOmega_n^2)\{1+\lito(1)\}=2\{n\tau(1-\tau)\}^{-2}\tr(\bOmega^2)\{1+\lito(1)\}$.
	In addition, under the given conditions, we have   $\hbeta_0/\beta_0\convP 1, \hbeta_1/\beta_1\convP 1$, and $\hd/d\convP 1$ as $n, p\to\infty$.
	We first prove (\ref{powerA.sec2}). Under Conditions C1--C3, Theorems~\ref{Tnexp.thm} and \ref{Tnlim.thm}   indicate that  as $n,p\to\infty$, we have  $\tT_{n,p,0}=T_{n,p,0}/\sigma_T\convL \zeta$ where $\zeta$ is defined in Theorem \ref{Tnlim.thm}. It follows that as $n, p\to\infty$,  we have
	\begin{equation}\label{power1.app}
		\begin{split}
			&\quad\Pr\left\{T_{n,p}\ge\hbeta_0+\hbeta_1\chi_{\hd}^2(\alpha)\right\}\\
			&= \Pr\left[\tT_{n,p,0} \ge \frac{\beta_0+ \beta_1 \chi^2_{d}(\alpha)}{\sigma_T}-\frac{\tr\left\{(\bSigma_1-\bSigma_2)^2\right\}}{\sqrt{2\{n\tau(1-\tau)\}^{-2}\tr(\bOmega^2)}}\right]\{1+\lito(1)\}  \\
			&=  \Pr\left[\zeta \ge \frac{\chi_{d}^2(\alpha)-d}{\sqrt{2d}}-\frac{n\tau(1-\tau)}{\sqrt{2\tr(\bOmega^2)}} \tr\left\{(\bSigma_1-\bSigma_2)^2\right\} \right]\{1+\lito(1)\},
		\end{split}
	\end{equation}
	where  $\tau$ is defined in Condition C1.
	
	We next prove (\ref{powerB.sec2}). Under the given conditions, when $d\to\infty$,  Theorem~\ref{Tnlim.thm}  indicates that  as $n,p\to\infty$, we have  $\tT_{n,p,0} \convL \zeta\sim N(0,1)$
	and  $\tT_{n,p,0}^* \convL N(0,1)$. In addition, as  $d\to\infty$, we have  $\{\chi_d^2(\alpha)-d\}/\sqrt{2d}\to z_{\alpha}$ where $z_{\alpha}$ denotes the upper $100\alpha$-percentile of $N(0,1)$. Then by (\ref{power1.app}), as $n,p\to\infty$,  we have
	\[
	\Pr\left\{T_{n,p}\ge\hbeta_0+\hbeta_1\chi_{\hd}^2(\alpha)\right\}
	=\Phi\left[-z_{\alpha}+\frac{n\tau(1-\tau)}{\sqrt{2\tr(\bOmega^2)}} \tr\left\{(\bSigma_1-\bSigma_2)^2\right\} \right]\{1+\lito(1)\},
	\]
	where $\Phi(\cdot)$ denotes the cumulative distribution of $N(0,1)$. The proof is complete.

		\bibliographystyle{apalike}
		\bibliography{HDmeanTest}
	\end{document}